\documentclass[11pt]{amsart}
\usepackage{amsfonts,amssymb,amsthm,amsmath,amsxtra,amscd,verbatim,eucal}
\usepackage[all]{xy}
\usepackage[dvips]{graphics}
\usepackage[usenames]{color}

\theoremstyle{plain}
\newtheorem{thm}{Theorem}[section]
\newtheorem*{thmA}{Theorem A}
\newtheorem*{thmB}{Theorem B}

\newtheorem{lem}[thm]{Lemma}
\newtheorem{prop}[thm]{Proposition}
\newtheorem{cor}[thm]{Corollary}

\theoremstyle{definition}
\newtheorem{defi}[thm]{Definition}

\newtheorem{eg}[thm]{Example}

\theoremstyle{remark}
\newtheorem{rmk}[thm]{Remark}
\newtheorem{notation}[thm]{Notation}

\newtheorem{step}{Step}

\numberwithin{equation}{section}

\def\Z{{\mathbb Z}}
\def\N{{\mathbb N}}

\def\C{{\mathbb C}}
\def\A{{\mathbb A}}
\def\R{{\mathbb R}}
\def\Q{{\mathbb Q}}
\def\P{{\mathbb P}}

\def\H{\mathcal{H}}

\def\O{\mathcal{O}}

\def\J{\mathcal{J}}

\def\W{\mathcal{W}}
\def\ZZ{\mathcal{Z}}

\def\aa{\mathfrak{a}}

\def\mm{\mathfrak{m}}

\def\xx{\underline{x}}%{{\bf x}}
\def\yy{\underline{y}}%{{\bf y}}

\def\a{\alpha}

\def\g{\gamma}
\def\d{\delta}
\def\f{\phi}
\def\ff{\psi}

\def\ep{\epsilon}

\def\l{\lambda}
\def\n{\nu}
\def\m{\mu}
\def\om{\omega}
\def\p{\pi}
\def\r{\rho}
\def\s{\sigma}
\def\t{\tau}

\def\x{\xi}
\def\z{\zeta}
\def\D{\Delta}
\def\G{\Gamma}
\def\L{\Lambda}
\def\S{\Sigma}

\def\.{\cdot}
\def\^{\widehat}
\def\~{\widetilde}
\def\o{\circ}
\def\ov{\overline}
\def\rat{\dashrightarrow}
\def\surj{\twoheadrightarrow}

\def\lru{\lceil}
\def\rru{\rceil}

\def\({\left(}
\def\){\right)}

\renewcommand{\subseteq}{\subset}

\renewcommand{\and}{ \ \ \text{ and } \ \ }
\renewcommand{\for}{ \ \ \text{ for } \ \ }
\newcommand{\fall}{ \ \ \text{ for all } \ \ }
\newcommand{\fevery}{ \ \ \text{ for every } \ \ }

\newcommand{\HH}[2]{H^{{#1}} \big({#2}   \big) }

\newcommand{\MI}[1]{\mathcal{J} \( {#1} \) }
\newcommand{\discr}[2]{a_{{#1}}\({#2}\)}

\newcommand{\can}[1]{{\rm can}\({#1}\)}

\DeclareMathOperator{\codim} {codim}
\DeclareMathOperator{\Hom} {Hom}
\DeclareMathOperator{\Spec} {Spec}
\DeclareMathOperator{\Proj} {Proj}
\DeclareMathOperator{\Aut} {Aut}
\DeclareMathOperator{\Bir} {Bir}

\def\Bl{{\rm Bl}}
\def\T{T}

\DeclareMathOperator{\ord} {ord}

\DeclareMathOperator{\val} {val}
\DeclareMathOperator{\Cont} {Cont}

\DeclareMathOperator{\ini} {in}
\DeclareMathOperator{\mld} {mld}
\DeclareMathOperator{\length} {length}
\DeclareMathOperator{\depth} {depth}

\begin{document}

%\vskip-1.5cm
%\begin{footnotesize}
%\hfill \jobname
%\end{footnotesize}
%\vskip1.5cm

\title{Birationally rigid hypersurfaces}

\author{Tommaso de Fernex}

\address{Department of Mathematics, University of Utah, 155 South 1400
East, Salt Lake City, UT 48112-0090, USA}
\email{{\tt defernex@math.utah.edu}}

\begin{abstract}
We prove that for $N \ge 4$, all smooth hypersurfaces of degree $N$
in $\P^N$ are birationally superrigid. 
First discovered in the case $N=4$ by Iskovskikh and
Manin in a work that started this whole direction of research,
this property was later conjectured to hold in general by Pukhlikov.
The proof relies on the method of maximal singularities in combination
with a formula on restrictions of multiplier ideals.
\end{abstract}


\thanks{The research was partially supported by NSF grants DMS-0111298, DMS-0548325
and CAREER grant DMS-0847059.}

\thanks{Compiled on \today. Filename {\tt \jobname}}

\maketitle

\section{Introduction}

Mori fiber spaces, alongside 
the minimal models, represent the terminal objects in the minimal 
model program. While minimal models can only be connected by chains
of flops, the birational pliability between Mori fiber spaces is in general
more rich. The notion of birational rigidity was introduced
to distinguish those Mori fiber spaces whose birational 
geometry is as simple as possible. 

Fano manifolds of Picard number one are,
in a trivial way, examples of Mori fiber spaces.
In this paper we address the question
asking which smooth Fano projective hypersufaces are
birationally rigid.
The complete answer to this question is given by the following result,
which was conjectured by Pukhlikov in \cite{Pu2}.

\begin{thmA}\label{intro-thm:bir-rigidity}
For $N \ge 4$, every smooth complex hypersurface
$X \subset \P^N$ of degree $N$ is birationally superrigid.
\end{thmA}

This results implies in particular that any such 
$X$ is nonrational, is the only Mori fiber space
in its birational equivalence class, and $\Bir(X) = \Aut(X)$.
According to a recent theorem of Odaka and Okada \cite{OO}, it also
implies that $X$ is slope stable in the sense of Ross and Thomas \cite{RT}.

Theorem~A has roots going back to the
work of Fano on the birational geometry of threefolds \cite{Fa1,Fa2}.
The case $N=4$ of Theorem~A is due to
Iskovskikh and Manin \cite{IM} for their proof of nonrationality
of smooth quartic threefolds in $\P^4$,
which, together with the theorem of 
Clemens and Griffiths on the nonrationality of smooth cubic threefolds
\cite{CG}, gave in the early seventies the first counter-examples to the L\"{u}roth
problem. Notably, the statement of Iskovskikh and Manin's theorem
already appears in Fano's work.
The powerful method introduced in \cite{IM},
known as the {\it method of maximal singularities},
uses techniques derived from the theory of resolution of singularities
to correct Fano's approach and hence successfully generalizes
Noether's inequality \cite{Noe} to quartic threefolds.

Considering the problem of birational rigidity for
a hypersurface $X \subset \P^N$, 
the question is only relevant when the degree of $X$ is $N$ and $N \ge 4$.
Thus Theorem~A characterizes all smooth complex projective hypersurfaces
that are birationally rigid.
Even though it was only explicitly conjectured in the late nineties,
this result was already suggested to be true---at least in some
weaker form---in the years
following the publication of \cite{IM} (cf. \cite{Pu2} for a discussion).

Using suitable extensions of the methods
of \cite{IM}, a few more low dimensional cases of Theorem~A were
established in a series of papers \cite{Pu1,Che,dFEM2},
and a proof of birational rigidity for
general hypersurfaces was given in \cite{Pu2}.
This result was then extended to include general hypersurfaces with generic
isolated singularities (a case which is not covered by Theorem~A) in \cite{Pu5}.
A related result was obtained in \cite{Ko1}, where it is proven that,
for $2 \lru (N+2)/3 \rru \le d \le N$, very
general hypersurfaces of degree $d$ in $\P^N$ are nonrational,
and in fact they are not birationally equivalent to any conic bundle
if $3 \lru (N+2)/4 \rru \le d \le N$.
An argument to prove the result stated in Theorem~A was proposed in \cite{Pu3},
but it contains a gap (cf. \cite[Remark~4.2]{dFEM2});
the same result is also claimed in \cite{Pu4},
but the given proof only applies to general hypersurfaces.

\medskip

The proof of Theorem~A uses the method of maximal
singularities---as any proof of birational rigidity given so far.
Following this method, one reduces to study log-pair structures on $X$ whose
singularities are strictly canonical, the delicate case being the one in which the
maximal singularity of the pair is centered at a point $P$ of $X$.

The main difficulty in higher dimensions is that base loci of
birational maps typically have large dimension.
At least in principle, one would like to cut down $X$ around $P$ until
the base locus becomes zero dimensional, as
in Corti's proof of Iskovskikh and Manin's theorem where
one restricts to a hyperplane section of the threefold \cite{Co2}.
The delicate part is to
keep track of how the singularities of the pair change in the process of cutting down.

One way one can try to do this is by taking into account the associated multiplier ideal.
The Connectedness Theorem of Shokurov and Koll\'ar implies that the singularities of
the pair are no longer log-terminal after
restricting to one hyperplane section: this corresponds to the appearance of a nontrivial
multiplier ideal. However, the same techniques do not apply
when further restrictions are taken because
multiplier ideals may ``stabilize'' after the first step,
as Example~\ref{eg:1} illustrates.
We suspect that this is the reason because this approach has not been followed
before.

To deal with this difficulty we make a careful selection of the hyperplane sections
and use arc spaces in place of vanishing theorems
to investigate how singularities change throughout the process.
The key technical result of the paper consists of an {\it inversion of adjunction} formula
on log-discrepancies of non-effective pairs, and is here presented
using the language of multiplier ideals.

Consider the affine space $\A^n$ with affine coordinates centered at a point $P$.
Let $Z \subset \A^n$ be a proper closed subscheme, and
suppose that, for some $c > 0$, the point $P$ is a center of non
log-terminality for the pair $(\A^n,cZ)$.
This means that there is a smooth variety $V$, a
proper birational morphism $f \colon V \to \A^n$, and a
prime divisor $E$ on $V$ such that $f(E) = P$ and
$$
c \val_E(I_Z) \ge \ord_E(K_{V/\A^n}) + 1,
$$
where $I_Z \subset \O_X$ is the ideal sheaf of $Z$ and $K_{V/\A^n}$
is the relative canonical divisor.
We can assume that $f$ factors through the blowup of $\A^n$ at $P$.
The center of $E$ on the blowup is contained in the exceptional
divisor, and thus its points determine tangent directions at $P$. 
Among these, most originate from arcs in $V$ that have contact with $E$;
these are what we call the {\it principal tangent directions of $E$ at $P$}
(cf. Proposition~\ref{lem:C_E-val_E(P)} and Definition~\ref{defi:principal-tangents}).
Using this terminology, we have the following result.

\begin{thmB}\label{intro-thm:discr-restr}
With the above notation, suppose that the valuation $\val_E$ is
invariant under the homogeneous $\C^*$-action on the coordinates of $\A^n$. 
Let $Y \subset \A^n$ be a positive dimensional linear subspace that is
not contained in $Z$ and is tangent to a 
principal tangent direction of $E$ at $P$.
Then, denoting by $\mm$ the maximal ideal of $P$ in $\A^n$
and by $\J(Y,cZ|_Y)$ the multiplier ideal of $(Y,cZ|_Y)$, we have
$$
\mm^{\codim(Y,\A^n)} \not\subseteq \J(Y,cZ|_Y).
$$
\end{thmB}

Membership to multiplier ideals can be expressed by
conditions on discrepancies of non-effective pairs.
It is precisely because of such non-effectiveness that the usual
methods relying on vanishing theorems cannot be employed. The proof of Theorem~B
uses the characterization of divisorial valuations
and log-discrepancies via contact loci in arc spaces given in \cite{ELM}.

In order to apply the above result to prove Theorem~A, we use generic projections
in a similar fashion as in \cite{dFEM2}
in combination with flat degenerations to homogeneous ideals,
where singularities are detected by a homogeneous valuations.
The most delicate part is to trace the
valuations throughout the degenerations. 
In particular, the necessary tangency condition on $Y$
imposed in Theorem~B will require us to investigate the effects that both
projections and degenerations have on the contact loci in arc spaces
that are associated to the valuations. 

Arc spaces have already been proved to be a powerful tool
in related problems (see, for example, the papers \cite{EMY,EM}
on inversion of adjunction).
A correspondence between exceptional divisors in resolution
of singularities and certain loci on arc spaces
was first explored by Nash \cite{Na}. The relationship between
arc spaces and singularities of pairs, implicit in the 
change-of-variable formula in motivic integration \cite{Ko,DL,Lo},
was later pointed out and investigated in the work of Musta\c{t}\v{a} \cite{Mu1}.
These correspondences have important implications to the birational geometry 
of algebraic varieties. 
The present paper gives a new application to a problem in classical algebraic geometry.

\medskip

The first part of the paper is devoted to the development of the necessary tools
from the geometry of arc spaces and contact loci leading to the proof of
Theorem~B.
The main result of the paper, Theorem~A, is then addressed in Section~6
where the question on the birational rigidity of Fano hypersurfaces is
reduced, via the Noether--Fano inequality, to a suitable
lower-bound to canonical thresholds (see Theorem~\ref{thm:bound-codim-2-case}).
The last two sections of the paper are devoted to the proof of such lower-bound.

\medskip

We work over the field of complex numbers.
With the term {\it variety} we will always 
refer to an integral separated scheme of finite type over $\Spec\C$, and
a {\it point} is intended to be a closed point unless stated otherwise.

\subsection*{Acknowledgement}

I would like to thank Roi Docampo and Mircea Musta\c{t}\v{a} for many 
enlightening discussions on some of the technical points of this paper. 
I am especially grateful to Musta\c{t}\v{a}
for pointing out a mistake in a previous version of this work. 
The paper has also benefited from conversations with Nero Budur
and Lawrence Ein.
I would like to thank the referee for a careful reading and
many valuable comments and corrections.

\section{Pairs and singularities}
\label{sect:1}

The theory of singularities of pairs and multiplier ideals plays a key role
in the proof of Theorem~A.
Standard references for these topics are \cite{Kol,KM,Laz}.
As we will adopt a slight different setting with respect to those
treated in these references, 
we briefly recall the basic definitions and terminology
as they will be used in the paper.

First of all, we will restrict ourselves to a setting 
where the ambient space is smooth. Let therefore $X$ be a smooth variety.
A divisor {\it over} $X$ is
a divisor $E$ on some smooth variety $X'$ endowed with a
proper birational morphism $f \colon X' \to X$. 
The divisor $E$ is said to be {\it exceptional} over $X$
if the image of its support has codimension greater than or equal to two.
When $E$ is a prime divisor, 
the image $f(E) \subset X$ is called the {\it center} of $E$ in $X$, 
and the {\it discrepancy} of $E$ with respect to $X$ is the number
$k_E(X) := \ord_E(K_{X'/X})$
where $K_{X'/X}$ is the relative canonical divisor of $f$
(since $X$ is smooth, $K_{X'/X}$ is locally defined by the Jacobian of $f$). 

A prime divisor $E$ over $X$ determines a {\it valuation} $\val_E$ on the function field  
of $X$. If $E$ is on $X'$ and $f\colon X' \to X$ is as above, this valuation
is defined by $\val_E(h) := \ord_E(h\o f)$ for every element $h \in \C(X)^*$.
Note that the center in $X$ of $\val_E$ is equal to the generic point of the center of $E$.
In general, if $\n$ is a valuation on the function field of $X$ having center in $X$
(e.g., $\n = \val_E$ for some divisor $E$ over $X$), $Z \subset X$
is a proper closed subscheme, and $I_Z\subset \O_X$ is the ideal sheaf of $Z$, 
then we denote 
$$
\n(Z) = \n(I_Z) := \min \{ \, \n(h) \mid 
\text{$h \in I_Z(U)$, $U \subset X$ open containing the center of $\n$ }\}.
$$

We consider pairs $(X,Z)$ consisting of a smooth variety $X$ and
a finite formal linear combination $Z = \sum c_iZ_i$ of
proper closed subschemes $Z_i$ of $X$ with real coefficients $c_i$.
This is where we diverge from the usual setting, where typically one
considers pairs of the form $(X,D)$ where $X$ is a normal variety and $D$ is a
$\Q$-divisor such that $K_X + D$ is $\Q$-Cartier ($K_X$ denoting a
canonical divisor on $X$). This last condition is used to define the appropriate
notion of discrepancy (see the above cited referenced for the precise definitions).
A general setting for singularities of pairs an multiplier ideals on normal varieties
which extends the usual setting and incorporates the one considered here is
developed in details in \cite{dFH}. 

As we always take $X$ to be a smooth variety, 
$K_X$ is already Cartier and we do not need to worry about the definition of discrepancy. 
On the other hand, working with subschemes (and formal linear combinations of these)
in place of divisors is very natural in order to study singularities
of rational maps, as these can be measured in terms of the singularities
of the base schemes of the linear systems defining the maps.
This is precisely the aim of this paper, and thus motivates us to
work in this setting. At least when the coefficients $c_i$
are rational, there is a connection between one setting (formal $\Q$-linear combinations of
closed subschemes) and the other ($\Q$-divisors on $X$). The connection is similar to the one
explained in \cite{Laz} in the case of 
a formal rational power of an ideal sheaf, which correspond in our language
to the case where $Z = c_1 Z_1$ is the rational multiple of a single
closed subscheme. 

Let $(X,Z)$ be a pair, with $X$ smooth and $Z = \sum c_iZ_i$ as above.
For a prime divisor $E$ over $X$, we denote by
$$
a_E(X,Z):= k_E(X) + 1 - \sum c_i\val_E(Z_i)
$$
the {\it log-discrepancy} of $E$ with respect to $(X,Z)$.
If $T \subseteq X$ is a closed subset, 
the infimum over all log-discrepancies $a_E(X,Z)$ where $E$ ranges among all prime 
divisors with center in $T$ is called the {\it minimal log-discrepancy}
of $(X,Z)$ in $T$, and is denoted by $\mld(T;X,Z)$. 
The pair $(X,Z)$ is said to be {\it terminal} if $a_E(X,Z) > 1$
for every exceptional prime divisor $E$ over $X$, and
{\it Kawamata log-terminal} if $a_E(X,Z) > 0$
for every prime divisor $E$ over $X$.
We say that $(X,Z)$ is terminal (resp. Kawamata log-terminal) {\it in dimension $k$}
if the corresponding condition on discrepancies holds for divisors $E$
with center in $X$ of dimension $\ge k$.
If all $c_i \ge 0$, then the {\it canonical threshold} of $(X,Z)$ is
the supremum of the values $c \in \R$ for which $(X,cZ)$ is terminal, and
is denoted by $\can{X,Z}$.

A proper birational morphism $f \colon X' \to X$
is said to be a {\it log-resolution} of $(X,Z)$ if
$X'$ is smooth, the exceptional locus of $f$ is a divisor, 
each $f^{-1}(Z_i)$ is a divisor on $X'$,
and the union of all these divisors
has simple normal crossing support.
Assuming that $c_i \ge 0$ for all $i$, given a log-resolution
$f \colon X' \to X$ of $(X,Z)$,
the {\it multiplier ideal} of $(X,Z)$ is the ideal sheaf
$$
\J(X,Z) :=
f_*\O_{X'}\left(\Big\lru K_{X'/X} - \sum c_i\.f^{-1}(Z_i)
\Big\rru \right) \subseteq \O_X,
$$
where for an $\R$-divisor $D$ on a smooth variety
we denote by $\lru D \rru$ its round-up (computed componentwise).
It is easy to see that the co-support of $\J(X,Z)$ is equal to the union
of the centers of all prime divisors $E$ over $X$
such that $a_E(X,Z) \le 0$; this set is also called
the {\it locus of non log-terminal singularities} of the pair $(X,Z)$.
Exactly as in the usual setting, $\J(X,Z)$ gives a scheme
structure to this locus.

The basic properties of multiplier ideals in the setting treated in \cite{Laz}
also hold in the setting considered here.
The same arguments can be applied in this setting; 
alternatively, one can refer to \cite{dFH}.

\section{Arc spaces and valuations}

In this section we review the characterization of
divisorial valuations via arc spaces due to \cite{ELM}.
We start by briefly recalling the basic notions
in the theory of arc spaces, referring to \cite{Mu1}
for more details.

For $m \ge 0$, the {\it $m$-th jet scheme} $J_m X$ of a scheme
of finite type $X$ is characterized by
$$
\Hom(\Spec A , J_m X) = \Hom(\Spec A[t]/(t^{m+1}),X)
$$
for every $\C$-algebra $A$.
The truncation maps $A[t]/(t^{m+1}) \to A[t]/(t^{m})$
give rise to affine morphisms $J_{m}X \to J_{m-1}X$.
The {\it arc space} $J_\infty X$ of $X$ is then defined as the inverse limit
of the projective system $\{J_{m}X \to J_{m-1}X\}_{m \ge 1}$. 
For $m \ge p \ge 0$ there are natural projections
$$
\p \colon J_\infty X \to X,
\quad
\p_m \colon J_\infty X \to J_mX,
\and
\p_{m,p} \colon J_mX \to J_pX
$$
given by truncation, and for every $m\in \N \cup \{\infty\}$
there is a (trivial) section $\iota_m\colon X \to J_mX$ that
at each point $P$ associates the constant jet or arc of $X$ at $P$.
For every subset $W \subseteq J_\infty X$
we denote $W_m := \p_m(W) \subseteq J_mX$.
A morphism of schemes $f \colon X \to Y$ induces, by composition, 
natural maps $f_m \colon J_mX \to J_mY$ that commute with the projections. 

For every $\l \in \C^*$, the automorphisms of the rings $A[t]/(t^{m+1})$
given by $t \mapsto \l t$ determine automorphisms of the schemes $J_mX$
that are compatible with the various projections and trivial sections.
The induced action on $J_\infty X$ will be denoted by $\g(t) \mapsto \g(\l t)$, where
$\g(t)$ stands for the typical arc.
Since we will consider other actions on jet schemes and arc spaces, 
we will refer to this action as the {\it canonical $\C^*$-action}.

\begin{rmk}\label{rmk:natural-action}
If $X = \A^n$ and $\xx = (x_1,\dots,x_n)$ are affine coordinates,
then $J_m\A^n$ is an affine space with
coordinates $(\xx,\xx',\dots,\xx^{(m)})$,
where $\xx^{(p)} = (x_1^{(p)},\dots,x_n^{(p)})$, and
the canonical $\C^*$-action of $J_\infty\A^n$
is given is these coordinates by $x_i^{(p)} \mapsto \l^p x_i^{(p)}$.
In general, if $X$ is an arbitrary smooth variety, 
one obtains a similar description of the fibers of $J_mX \to X$
by fixing a system of parameters $\xx = (x_1,\dots,x_n)$
in the maximal ideal of a point $P$ in $X$: the fiber over $P$
is an affine space with
coordinates $(\xx',\dots,\xx^{(m)})$
and the action on this fiber is again given by $x_i^{(p)} \mapsto \l^p x_i^{(p)}$.
\end{rmk}

In the remainder of this section, we assume that $X$ is
a smooth variety.

\begin{defi}
A {\it cylinder} of $J_\infty X$ is
a subset $C \subseteq J_\infty X$ that is the inverse image
of a constructable set on some finite level $J_pX$.
\end{defi}

\begin{rmk}\label{rmk:C_m-constr-irred}
By a theorem of Chevalley, for every cylinder $C \subseteq J_\infty X$ and
every $m \ge 0$ the subset $C_m$ of $J_m X$ is constructable.
Moreover, if $C$ is irreducible, then the closure of $C_m$ in $J_m X$
is irreducible.
\end{rmk}

Consider now a proper closed subscheme $Z \subset X$.

\begin{defi}
The {\it order of contact} of an arc $\g \in J_\infty X$ along $Z$ is defined as
$$
\ord_\g(Z) = \ord_\g^X(Z)
:= \sup \big\{\, m  \mid \g^*(I_Z) \subseteq (t^m)\C[[t]] \,\big\}
\in \Z \cup \{\infty\},
$$
where $I_Z \subseteq \O_{\p(\g),X}$ is the ideal locally defining $Z$ near $\p(\g)$.
\end{defi}

\begin{rmk}\label{rmk:order-restriction}
If $Y \subseteq X$ is a submanifold that is not contained in $Z$, and
$\g \in J_\infty Y \subseteq J_\infty X$, then $\ord_\g^Y(Z|_Y) = \ord_\g^X(Z)$.
\end{rmk}

For every integer $q \ge 0$, let
$$
\Cont^q(Z) = \{\, \g \in J_\infty X \mid \ord_\g(Z) = q \,\},
\quad
\Cont^{\ge q}(Z) = \{\, \g \in J_\infty X \mid \ord_\g(Z) \ge q \,\}.
$$
Each of these sets is a cylinder in $J_\infty X$.

\begin{defi}
Any subset $W \subset J_\infty X$ that is
the closure of an irreducible component of $\Cont^q(Z)$ for some
proper closed subscheme $Z \subset X$ and some $q \ge 1$
is said to be an {\it irreducible contact locus}.
\end{defi}

\begin{rmk}
Every irreducible contact locus is an irreducible closed cylinder of $J_\infty X$,
does not dominate $X$, and is invariant under the canonical $\C^*$-action.
\end{rmk}

\begin{defi}
A {\it divisorial valuation} of $\C(X)$ is a valuation of the
form $q\val_E$ for some positive integer $q$ and some prime divisor $E$ over $X$.
\end{defi}

\begin{defi}
Associated to any non-empty closed and irreducible cylinder $C \subset J_\infty X$
that does not dominate $X$, there is a discrete valuation $\val_C$ of $\C(X)$
defined by the condition $\val_C(h) = \min\{\ord_\g(h) \mid \g \in C\}$
for any element $h \in \C(X)^*$.
Every valuation constructed in this way is called a {\it cylinder valuation}.
If the cylinder is equal to an irreducible
contact locus $W$ in $J_\infty X$, then we say that $\val_W$ is a {\it contact valuation}.
\end{defi}

\begin{rmk}
With the same notation as in the last definition,
we have $\val_C(h) = \ord_\g(h)$ for a general $\g \in C$
(the generality may depend on $h$).
\end{rmk}

One of the main results in \cite{ELM} says that
the set of divisorial valuations, the set of cylinder valuations,
and the set of contact valuations are all the same subset
of the set of all valuations of the function field of $X$.
More precisely, starting with a divisorial valuation
of $\C(X)$, we fix a proper birational map
$f \colon X' \to X$, where $X'$ is smooth and contains a 
smooth prime divisor $E$ such that
the valuation is written as $q\val_E$ for some $q \ge 1$.
We define
$$
W^q(E) := \ov{f_\infty\big(\Cont^q(E)\big)} \subset J_\infty X.
$$
It is easy to see that the set $W^q(E)$ is irreducible
and only depends on the valuation $q\val_E$.
The assertions in the next statement that are not explicitly stated
in \cite[Theorems~A and~C]{ELM} can be found in the proofs therein.

\begin{thm}[\protect\cite{ELM}]\label{thm:ELM}
With the above notation, the following properties hold. 
\begin{enumerate}
\item[(i)]
$W^q(E)$ is an irreducible cylinder of codimension
$$
\codim(W^q(E),J_\infty X) = q(k_E(X) + 1),
$$
and $\val_{W^q(E)} = q\val_E$.
\item[(ii)]
Every irreducible contact locus $W$ is equal to $W^q(E)$ for some 
prime divisor $E$ over $X$ and some integer $q\ge 1$.
\item[(iii)]
Every cylinder valuation $\val_C$ is equal to a contact
valuation $\val_W$ for some irreducible contact locus $W$ containing $C$. In particular,
if $E$ and $q$ are such that $\val_C = q\val_E$, then
$$
\codim(C,J_\infty X) \ge q(k_E(X) + 1).
$$
\end{enumerate}
\end{thm}

\begin{rmk}
It follows from the theorem that each $W^q(E)$ is an irreducible contact locus.
\end{rmk}

\section{Higher order tangent spaces and principal tangent directions}

Let $X$ be a smooth variety of dimension $n \ge 1$, and fix a point $P \in X$.
For any valuation $\n$ on the function field of $X$ 
we denote by $\n(P) := \n(\mm_P)$ the valuation of the 
maximal ideal $\mm_P \subset \O_X$ of $P$
(this is consistent with the notation introduced in
Section~\ref{sect:1} if $P$ is regarded as a subscheme of $X$).
For every $m \ge 0$, let $P_m \in J_mX$ denote the constant $m$-th jet
of $X$ at $P$.

\begin{defi}\label{defi:T^(m)X}
For $m \ge 1$, we define the {\it $m$-th order tangent space} of $X$ at $P$
and the {\it $m$-th order tangent bundle} of $X$ to be, respectively,
$$
\T_P^{(m)}X := \p_{m,m-1}^{-1}\big(P_{m-1}\big)
\and
\T^{(m)}X := \p_{m,m-1}^{-1}\big(\iota_{m-1}(X)\big).
$$
\end{defi}

We have $\T'X = \T X$, the tangent bundle of $X$.
In general, the scheme $\T^{(m)}X$ is characterized by
$$
\Hom(\Spec A , \T^{(m)}X) = \Hom(\Spec A[t^m]/(t^{m+1}),X)
$$
for every $\C$-algebra $A$. The embedding in $J_mX$
corresponds to the inclusion $A[t^m]/(t^{m+1}) \subseteq A[t]/(t^{m+1})$.
If $m,p \ge 1$, then the ring isomorphism
$A[t^m]/(t^{m+1}) \to A[t^p]/(t^{p+1})$ mapping $a + bt^m$ to $a + bt^p$
($a,b \in A$) gives rise to a natural isomorphism
$$
\ff_{m,p} \colon \T^{(m)}X \to \T^{(p)}X.
$$
Clearly $\p_{p,0}\o\ff_{m,p} = \p_{m,0}$ and
$\ff_{m,p} = \ff_{q,p}\o\ff_{m,q}$ for every $m,p,q \ge 1$.
In particular, the vector bundle structure of $\T X$ naturally
lifts, via $\ff_{m,1}$, to a vector bundle structure on $\T^{(m)}X$,
each $\ff_{m,p}$ becomes an isomorphism of vector bundles,
and each $\T_P^{(m)}X$ a vector space.

\begin{lem}\label{lem:T_P^mX-ff-Y_m}
If $Y \subseteq X$ is a closed subscheme, then
$J_mY \cap \T_P^{(m)}X$ is a linear subspace of $\T_P^{(m)}X$ and
$\ff_{m,p}\big(J_mY \cap \T_P^{(m)}X\big) = J_pY \cap \T_P^{(p)}X$,
for every $m,p \ge 1$.
\end{lem}

\begin{proof}
Given a system of parameters $\xx = (x_1,\dots,x_n)$ of $P$ in $X$,
we have $\T^{(m)}_PX = \Spec \C[\xx^{(m)}]$,
the coordinates $\xx^{(m)}$ are linear on $\T_P^{(m)}X$,
and $\ff_{m,p}$ is given by the equations
$\xx^{(p)} = \xx^{(m)}$.
Observe that if $h(\xx)$ is a nonconstant homogeneous
polynomial, then for every $0 \le q < m$ its $q$-th differential
$h^{(q)}(\xx,\xx',\dots,\xx^{(q)}) :=
D^q(h(\xx))$
vanishes identically over $\T^{(m)}_PX$, and it also vanishes when $q=m$
if $\deg(h) \ge 2$. On the other hand, if
$\ell(\xx)$ is a linear form, then
$\ell^{(m)}(\xx,\xx',\dots,\xx^{(m)}) =
\ell(\xx^{(m)})$ on $\T^{(m)}_PX$.

We can assume that $P \in Y$ and $Y \ne X$.
In the formal neighborhood of $X$ at $P$, $Y$ is defined by equations
$f_\a(\xx)=0$ where $f_\a \in \C[[\xx]]$, $1 \le \a \le r$.
Let $\ell_\a(\xx)$ denote the linear part of $f_\a(\xx)$.
Then $J_mY \cap \T_P^{(m)}X$ is defined in $\T_P^{(m)}X$ by the equations
$\ell_\a(\xx^{(m)}) = 0$.
Both assertions follow.
\end{proof}

\begin{lem}\label{lem:W_m-in-T^m}
Let $W \subset X_\infty$ be an irreducible contact locus.
Then for every $m \ge 1$ the set
$W_m \cap \T_P^{(m)}X$ is a homogeneous subset of $\T_P^{(m)}X$.
\end{lem}

\begin{proof}
Take an arbitrary point $\g(t) \in W_m \cap \T_P^{(m)}X$.
Since $W$ is invariant under the canonical $\C^*$-action, 
$\g(\l t) \in W_m$ for every $\l \in \C^*$. The lemma follows by observing that
$\g(\l t) = \l^m\.\g(t)$ in $\T_P^{(m)}X$.
\end{proof}

For any closed set $C \subset J_\infty X$ 
that is contained in the fiber $\p^{-1}(P)$ over $P$, we denote 
$$
\m_P(C) := \min \{\, m \ge 0 \mid C_m \ne P_m \,\}.
$$
Note that $\m_P(C)$, if finite, is equal to the unique integer $\m$ such that
$P_\m \ne C_\m \subseteq \T_P^{(\m)} X$.

\begin{lem}
\label{lem:3.4}
If $C$ is an irreducible closed cylinder contained in $\p^{-1}(P)$
and $\val_C$ is the associated valuation, then 
$$
\m_P(C) = \val_C(P).
$$
\end{lem}

\begin{proof}
Fix a local system of parameters $\underline{x}$ in the maximal ideal of $P$.
If $(a_i^{(p)})$ are the coordinates of a general arc $\g \in C$ in the 
coordinate system $(x_i^{(p)})$, then
$\m_P(C)$ is the smallest integer $\m$ such that $a_i^{(\m)} \ne 0$ for some $i$. 
On the other hand, if we taken a general linear form $\ell(\underline{x})$, 
we see that $\ord_\g(\ell(\underline{x}))$ is equal to the same number $\m$. 
\end{proof}

\begin{notation}\label{not:G_E}
Let $E$ be a prime divisor over $X$ with center $P$.
Assuming without loss of generality that $E$ is a divisor on 
a model dominating the blowup $\Bl_PX$ of $X$ at $P$,
we denote by $\G_E$ the center of $E$ in $\Bl_PX$. Since $\G_E$ is contained
in the exceptional divisor of the blowup, we can take its cone in $\T_PX$;
we denote this cone by $\^\G_E$.
\end{notation}

\begin{prop}\label{lem:C_E-val_E(P)}
Let $W\subseteq J_\infty X$ be an irreducible contact locus such that $\val_W = q\val_E$
for some prime divisor $E$ over $X$ with center $P$ and some integer $q \ge 1$, 
and let $\m = \m_P(W)$. Then $\ff_{\m,1}(W_\m) \subseteq \T_PX$
is a dense constructable subset of $\^\G_E$.
\end{prop}

\begin{proof}
The fact that $W_\m$ is a constructable subset of $T_P^{(\m)}X$ follows
by Chevalley's theorem, and implies that $\ff_{\m,1}(W_\m)$
is constructable in $\T_PX$. We need to show that $\ff_{\m,1}(W_\m)$
is contained in $\^\G_E$ and is dense in there.
First observe that $\^\G_E \subseteq \T_PX$ is defined by a homogeneous ideal,
and so is the closure of $W_\m$ in $\T_P^{(\m)}X$ (see Lemma~\ref{lem:W_m-in-T^m}).
Therefore, in order to check that $\ff_{\m,1}(W_\m)$ is a dense subset of $\^\G_E$,
it suffices to compare vanishing of homogeneous polynomials in
a local system of parameters $\underline{x}$ in the maximal ideal of $P$.
Fix an arbitrary homogeneous polynomial $h$ of degree $d \ge 1$
in $n$ variables. We first observe that, by Lemma~\ref{lem:3.4},
$$
q\val_E(h(\underline{x})) > \m d
\quad \Leftrightarrow \quad h\big(\underline{x}')|_{\^\G_E} \equiv 0.
$$
Now let $\g \in W$ be a general arc, so that in particular
$\val_W(h(\underline{x})) = \ord_\g(h(\underline{x}))$.
Note that $\p_\m(\g)$ is a general point of $W_\m$ and $\p_{\m-1}(\g) = P_{\m-1}$.
Then, if we denote by
$\underline{a}^{(\m)} = (a_1^{(\m)},\dots,a_n^{(\m)})$ the coordinates
of $\p_\m(\g)$ in $\T^{(\m)}_PX$ with respect to the system of coordinates
$\underline{x}^{(\m)}$, we see that
$$
\val_W(h(\underline{x})) > \m d \quad \Leftrightarrow \quad
h(\underline{a}^{(\m)}) = 0
\quad \Leftrightarrow \quad h(\underline{x}^{(\m)})|_{W_\m} \equiv 0.
$$
Since $q\val_E = \val_W$ and $\^\G_E$ is closed, 
we conclude that $\ff_{\m,1}(W_\m)$ is a dense subset of $\^\G_E$.
\end{proof}

Under the assumptions of the previous proposition, 
$\ff_{\m,1}(W_\m)$ is a homogeneous subset of $\T_PX$, 
and thus it is the cone over a constructable dense subset of $\G_E$.
We introduce the following terminology.

\begin{defi}\label{defi:principal-tangents}
Let $W = W^1(E) \subset J_\infty X$ where $E$ is a prime divisor over $X$ with center $P$, 
and let $\m := \m_P(W) = \val_E(P)$.
The non-zero elements in $\ff_{\m,1}(W_\m) \subseteq \^\G_E$ are said to be the 
{\it principal tangent vectors of $E$ at $P$}.
The elements in $\G_E$ given by homogeneous classes of 
non-zero elements in $\ff_{\m,1}(W_\m)$
are the {\it principal tangent directions of $E$ at $P$}. 
\end{defi}

\section{Homogeneous action and degeneration of contact loci}

Throughout this section, let $\A^n = \Spec\C[\xx]$
where $\underline{x} = (x_1,\dots,x_n)$, and
let $P \in \A^n$ denote the origin in these coordinates. 
We consider the homogeneous $\C^*$-action of $\A^n$ 
given by $x_i \mapsto s x_i$, for $s \in \C^*$. 
For short, we denote the action by $\xx \mapsto s\xx$.

The action induces a flat degeneration in $\C[\xx]$
to homogeneous forms by taking limit when $s \to 0$. 
Under this action, a polynomial $h$
degenerates to the homogeneous polynomial consisting of all the terms of
largest degree of $h$.
We denote by $\ini(\aa)$ the
corresponding flat degeneration of an ideal 
$\aa\subseteq \C[\xx]$.
If $Z \subset \A^n$ is a closed subscheme, we denote by $Z^0$
the homogeneous scheme defined by the ideal $\ini(I_Z)$. 
Once the action is clear from the context, 
we will refer to $Z^0$ as the {\it homogeneous limit} of $Z$. 

\begin{rmk}
\label{rmk:Z^0-cone}
If $\A^n = \P^n \smallsetminus H$ where $H \subset \P^n$ is a hyperplane, 
then a homogeneous $\C^*$-action is already determined 
by the choice of an origin $O \in \A^n$, 
and extends to a $\C^*$-action on $\P^n$ pointwise fixing $H$.
If $\ov Z \subset \P^n$ is the projective closure of a
closed subscheme $Z \subset \A^n$, then 
the homogeneous limit $Z^0$ always contains the affine
cone over $\ov Z \cap H$ with center in $O$. The two cones are in fact equal whenever
$\ov Z$ is a positive dimensional complete intersection in $\P^n$: indeed in this case 
$\ov Z \cap H$ is a complete intersection in $H$, cut out by forms
of the same degree as those cutting out $\ov Z$, and
one checks the equality between $Z^0$ and the cone over $\ov Z \cap H$ by
comparing their Hilbert polynomials. Note in particular that
in this case $Z^0$ is a complete intersection subscheme of $\A^n$. 
In general however, it is not true that $Z^0$ is the cone over $\ov Z \cap H$:
this fails if $Z$ is zero-dimensional or contains
zero-dimensional embedded points, since these are not detected by $\ov Z \cap H$.
\end{rmk}

The degeneration can be described geometrically, for an arbitrary closed subscheme
$Z \subset \A^n$, as follows. 
Let $\D = \Spec\C[s]$ and $\D^* = \D \smallsetminus\{0\}$. Let
$\ZZ^* \subseteq \C^*\times \A^n$ be the image of $\C^*\times Z$ under the map
$$
g\colon \D^* \times \A^n \to \D^*\times \A^n, \quad (s,\xx) \mapsto (s, s \xx),
$$
and denote by $\ZZ \subseteq \D \times \A^n$ the flat closure of $\ZZ^*$ 
under the inclusion $\D^* \times \A^n \subset \D \times \A^n$. 
Then $Z^0$ is simply given by the fiber $\ZZ_0$
of $\ZZ$ over $0 \in\D$ via the isomorphism
$\{0\} \times \A^n \to \A^n$ induced by the projection on the second factor. 
Note that $\ZZ_s \cong Z$ for every $s \in \D^*$.

This construction is compatible with general linear projections. 

\begin{prop}\label{prop:cone-proj}
Let $\f \colon \A^n \to \A^m$ be a surjection that is linear and homogeneous with respect to 
some given affine coordinates on $\A^n$ and $\A^m$. We consider the homogeneous
$\C^*$-actions associated to these coordinates. 
Let $Z \subset \A^n$ be a closed subscheme, and assume that the closure
$\ov{Z}$ of $Z$ in $\P^n$ is disjoint from the indeterminacy locus 
of the induced rational map $\ov{\f}\colon \P^n \rat \P^m$. Then
$$
\f(Z^0) = \f(Z)^0,
$$
where $Z^0 \subset \A^n$ and $\f(Z)^0 \subset \A^m$ are the homogeneous
limits of $Z$ and $\f(Z)$ as described above. 
\end{prop}

\begin{proof}
Let $\ZZ \subset \D \times \A^n$ be the subscheme giving the flat deformation of $Z$
to $Z^0$, and $\ZZ' \subset \D \times \A^m$ 
be the subscheme giving the flat deformation of $\f(Z)$
to $\f(Z)^0$. For every $s \ne 0$ the fiber $\ZZ_s$ maps to the corresponding
fiber $\ZZ'_s$.
By hypothesis, the closure $\ov{\ZZ}$ of $\ZZ$ inside $\D \times \P^n$ 
is disjoint from the indeterminacy locus of 
$$
1_\D \times \ov{\f} \colon \D \times \P^n \rat \D \times \P^m,
$$ 
and thus $(1_\D \times \ov{\f})(\ov{\ZZ})$ is a closed subscheme of $\D \times \P^m$. 
Restricting to $\D \times \A^m$, we conclude that 
$(1_\D \times \f)(\ZZ)$ is closed in $S \times \A^m$.
By construction, $\ZZ$ is flat over $\D$, that is, its ring of functions is 
a torsion-free $\C[s]$-module. This implies that the ring of functions
of $(1_\D \times \f)(\ZZ)$ is a torsion-free $\C[s]$-module, and hence 
$(1_\D\times \f)(\ZZ)$ is also flat over $\D$. 
Since $(1_\D \times \f)(\ZZ)$ agrees with $\ZZ'$ over $\D^*$, 
we conclude that $(1_\D \times \f)(\ZZ) = \ZZ'$. 
Restricting to the fiber over $s$, we see that $\f$ maps $Z^0$ to $\f(Z)^0$.
\end{proof}

The $\C^*$-action of $\A^n$ naturally lifts to $\C^*$-actions on 
$J_m\A^n = \Spec\C[\underline{x},\dots,\underline{x}^{(m)}]$ for each $m$,
given by $x_i^{(p)} \mapsto s x_i^{(p)}$, as well as on $J_\infty \A^n$. 
These actions are compatible with the various projections and trivial sections.
With slight abuse of notation,
once the coordinates $x_i^{(p)}$ are clearly specified
we will denote the action on $J_\infty \A^n$ by
$\g(t) \mapsto s \. \g(t)$. 
To distinguish it from the canonical $\C^*$-action $\g(t) \mapsto \g(\l t)$, 
will refer to it as the {\it homogeneous $\C^*$-action}. 

For every integer $m$, 
we obtain flat degenerations of closed subschemes of $J_m\A^n$ exactly as described
for subschemes of $\A^n$. 
We will use the same construction all the way up on $J_\infty \A^n$
to define the degeneration of a closed subset $W \subseteq J_\infty \A^n$
to its (set theoretic) {\it homogeneous limit} $W^0 \subseteq J_\infty \A^n$. 
The degeneration is obtained as follows. 

The map $g \colon \D^* \times \A^n \to \D^*\times \A^n$ defined above lifts to the morphism
$$
g_\infty \colon \D^*\times J_\infty \A^n \to \D^*\times J_\infty \A^n,
\quad (s,(x_i^{(p)})) \mapsto (s,s x_i^{(p)}).
$$
For any closed subset $W \subseteq J_\infty \A^n$, we consider the image
$\W^* = g_\infty(\D^*\times W) \subseteq \D^*\times J_\infty \A^n$ 
and let $\W \subseteq \D\times J_\infty \A^n$ be its closure . 
Then $W^0$ is given by the fiber $\W_0$ of $\W$ over $0 \in\D$. 
Note that if $W$ is invariant under the canonical $\C^*$-action then so is $W^0$, 
and similarly if $W$ is contained in the fiber over $P$ then so is $W^0$.

\begin{lem}\label{lem:W}
If $W \subseteq J_\infty \A^n$ is a nonempty closed cylinder and 
$m$ is such that $W = \p_m^{-1}(W_m)$ then, set-theoretically, $(W^0)_m = (W_m)^0$
and $W^0 = \p_m^{-1}((W_m)^0)$. In particular, $W^0 \ne \emptyset$.
\end{lem}

\begin{proof}
Both assertions follow by the fact that if $W$ is a cylinder over $W_m$,
then the equations cutting out $W$ in $J_\infty \A^n$
are polynomials in $(\ov{x},\dots, \ov{x}^{(m)})$. 
\end{proof}

\begin{prop}\label{prop:W}
With the above notation, if $W \subseteq J_\infty \A^n$ is
a nonempty closed irreducible cylinder in $J_\infty \A^n$, then
every irreducible component $C$ of $W^0$ has codimension
$$
\codim(C,J_\infty \A^n) = \codim(W,J_\infty \A^n).
$$
\end{prop}

\begin{proof}
Fix $m$ such that $W$ is a cylinder over $W_m$. 
By Lemma~\ref{lem:W}, $(W_m)^0$ is nonempty and
$W^0$ is a cylinder over $(W_m)^0$.
Therefore it suffices to look at the dimension of the irreducible components of $(W_m)^0$. 
Since $\W_m \to \D$ is surjective,  
every irreducible component $C_m$ of $(\W_m)^0$ has dimension 
$\dim C_m \ge \dim \W_m - 1$, and in fact equality holds since $\W_m$ is flat over $\D$. 
\end{proof}

Suppose now that $W = W^1(E) \subset J_\infty \A^n$
for some prime divisor $E$ over $\A^n$. 
We can look at the degeneration of $W$ into $W^0$
as providing a way of ``degenerating'' the valuation $\val_E$ to homogeneous valuations.
Specifically, we pick any irreducible component $C$ of $W^0$ and consider the associated 
valuation $\val_C$.
This is a homogeneous valuation since $W^0$ is invariant under the homogeneous
$\C^*$-action.
By Theorem~\ref{thm:ELM}, there is a prime divisor $F$ over $\A^n$
and an integer $q \ge 1$ such that $\val_C = q\val_F$. 
Note that if $E$ has a center $P$ then so does $F$ by 
Lemma~\ref{lem:W}.

\begin{prop}\label{prop:semicont}
With the above notation, 
let $Z$ be a proper subscheme of $\A^n$, fix $c > 0$, and consider the homogeneous
limit $Z^0\subset \A^n$ of $Z$. Then 
$$
q \,a_F(\A^n,cZ^0) \le a_E(\A^n,cZ).
$$
In particular, if $a_E(\A^n,cZ) \ge 0$, then $a_F(\A^n,cZ^0) \le a_E(\A^n,cZ)$.
\end{prop}

\begin{proof}
If $W$ and $C$ are as above, then we have
$$
q (k_F(\A^n) + 1) \le \codim(C,J_\infty\A^n) = \codim(W,J_\infty\A^n) = k_E(\A^n) + 1
$$
by Theorem~\ref{thm:ELM} and Proposition~\ref{prop:W}.
Let $\ZZ \subset \D \times \A^n$ be the flat family giving the deformation, 
so that $Z^0 = \ZZ_0$. Then 
$$
\val_E(Z) = \val_W(Z) = \val_{\W_s}(\ZZ_s) \le \val_C(Z^0) = q\val_F(Z^0)
$$
for every $s \ne 0$, by Remark~\ref{rmk:order-restriction},
semi-continuity, and the inclusion $C\subseteq \W_0 \cap J_\infty \A^n$.
Putting together, we have $q\, a_F(\A^n,cZ^0) \le a_E(\A^n,cZ)$.
The last assertion follows from the fact that $q \ge 1$.
\end{proof}

\begin{rmk}
It would be interesting to know more about the relationship between the valuations
$\val_E$ and $\val_F$ described in the above contruction. For instance, 
if $E$ has center $P$ in $\A^n$ 
and we take take the degeneration as $s \to \infty$,
which gives the degeneration to the normal cone, and $F$ is constructed 
in an analogous manner under this degeneration, 
then it is possible to prove that $E$ and $F$ have the same set of 
principal tangent directions at $P$. 
This however may fail for the degeneration given by $s \to 0$. 
A related question asks whether there is some 
relationship between $W^0$ and $(W')^0$ if $W' \subset J_\infty \A^n$ is any closed 
irreducible cylinder giving the same valuation $\val_{W'} = \val_W$.
Principal tangent directions will play a crucial role in the following sections. Potentially, 
a better control on these directions
could simplify some of the technical points in the proof 
of the main result of the paper. 
\end{rmk}

\section{Log-discrepancies and restrictions}

We begin by quickly reviewing the proof of the inversion of adjunction
formula in the smooth setting
using the arc interpretation of valuations and log-discrepancies
discussed in Section~2. We refer to \cite{EMY,EM,Kaw,EM2,Ish,dFD}
for the original formulation and for more general settings.

\begin{thm}[Inversion of Adjunction \protect{\cite{EMY,EM}}]\label{thm:usual-adj}
Let $X$ be a smooth complex variety, and let
$E$ be a prime divisor over $X$, with center $T \subset X$.
Let $Y \subset X$ be a smooth positive dimensional subvariety
of codimension $e$ containing $T$. 
Then there exist a divisor $F$ over $Y$ with center contained in $T$,
and an integer $q \ge 1$, such that
for any proper closed subscheme $Z \subset X$ not containing $Y$ and any $c \ge 0$
$$
q\,a_F(Y,cZ|_Y) \le a_E(X,cZ+eY).
$$
In particular, if $a_E(X,cZ+eY) \ge 0$, then $a_F(Y,cZ|_Y) \le a_E(X,cZ+eY)$.
\end{thm}

\begin{proof}
Let $W = W^1(E) \subset J_\infty X$ be the contact locus associated to $E$.
By Theorem~\ref{thm:ELM}, we have
$\codim(W,J_\infty X) = k_E(X)+1$.
Note that $J_\infty Y \cap W \ne \emptyset$ since they both contain the 
constant arcs through $T$.
Clearly $J_\infty Y \cap W$ is a cylinder in $J_\infty Y$;
let $C$ be any of its irreducible components.
Since at least $e\val_E(Y)$ of the equations defining $J_\infty Y$
in $J_\infty X$ already vanish identically on $W$, we have
$$
\codim(C,J_\infty Y) \le \codim(W,J_\infty X) - e\val_E(Y).
$$
By Theorem~\ref{thm:ELM}, there exist a prime divisor $F$
over $Y$ and a positive integer $q$ such that
$\val_C = q\val_F$ as valuations of $\C(Y)$, and
$\codim(C,J_\infty Y) \ge q(k_F(Y)+1)$. On the other hand, 
$q\val_F(Z|_Y) \ge \val_E(Z)$
(this assertion will be discussed in more details in the proof of
Theorem~\ref{thm:discr-restr:homog-case}).
Putting together, we obtain
\begin{align*}
q\,a_F(Y,cZ|_Y)
& = q\big(k_F(Y) + 1 - c\val_F(B|_Y)\big) \\
& \le k_E(X) + 1 - c\val_E(Z) - e\val_E(Y) \\
& = a_E(X,cZ + eY).
\end{align*}
To conclude, just observe that the center of $F$
in $Y$ is coincides with the image of $C$, and this is contained in $T$.
\end{proof}

\begin{rmk}
Inversion of adjunction is usually expressed
in terms of minimal log-discrepancies. In the setting of the theorem,
it states that
$$
\mld(T;Y,cZ) = \mld(T;X,cZ + eY)
$$
where $T \subset Y$ is any closed irreducible subset of codimension at least 2
and $Z \subset Y$ is a proper
closed subscheme. The two formulations are equivalent. 
The fact that the formula in Theorem~\ref{thm:usual-adj} follows
from this identity on minimal log-discrepancies is clear.
Conversely, while the inequality $\mld(T;Y,cZ) \ge \mld(T;X,cZ + eY)$ 
is an immediate consequence of the adjunction formula for the canonical class
applied along a common log-resolution, the reverse inequality can easily be deduced
by the statement of Theorem~\ref{thm:usual-adj} using the well-known fact that if
$\mld(T;Y,cZ) < 0$ then $\mld(T;Y,cZ) = -\infty$. 
The advantage of stating the theorem the way we did is that it allows us to focus on the
valuations computing (or bounding) the minimal log-discrepancy. 
\end{rmk}

Our goal in this section is to generalize the above theorem.
Roughly speaking, we would like to ``move'' the contribution of
the valuation of $Y$ from one side to the other in the formula.
This becomes an inversion of adjunction problem for non-effective pairs. 
Conditions on log-discrepancies can be translated in terms of 
multiplier ideals. From this point of view,
what we have in mind is a formula able to measure, in some
effective way, possible jumps on multiplier ideals through restrictions.

The next example brings to light some of the difficulties
hidden in this problem.

\begin{eg}\label{eg:1}
Let $X = \A^n$ with coordinates $(x_1,\dots,x_n)$. Fix pairwise coprime
positive integers $a_i$ ($n \ge i \ge 1$), starting with $a_n$ and then
increasing rapidly enough:
\begin{equation}\label{eq:eg:1}
0 < a_n \ll a_{n-1} \ll \dots \ll a_1.
\end{equation}
Then let $Z \subset X$ be the subscheme defined by the
integral closure of the ideal
$(x_1^{a_1},\dots,x_n^{a_n})$. Fix $n-1$ general linear forms on $X$,
and for every $e \in \{0,1,\dots,n-1\}$
let $Y_e \subseteq X$ be the linear subspace of codimension $e$
defined by the vanishing of the first $e$ of these forms.
After possibly a linear change of coordinates
(which would not change the given description of the ideal of $Z$),
we can assume that
$$
Y_e = \{ x_1 = \dots = x_e = 0\} \fevery 1 \le e \le n-1.
$$
Note that, in the induced coordinates $(x_{e+1},\dots,x_n)$ of $Y$,
$Z|_{Y_e}$ is defined by
the integral closure of the ideal $(x_{e+1}^{a_{e+1}},\dots,x_n^{a_n})$.
Let
$$
c_e = \frac{1}{a_{e+1}} + \dots + \frac{1}{a_n} \and c = c_0.
$$
There is a prime divisor $E$ over $X$ such that $\discr{E}{X,cZ} = 0$
(such divisor is determined by the toric blowup with weights
$a/a_i$ along the coordinates $x_i$, where $a := a_1\dots a_n$), and thus the multiplier ideal
$\J(X,cZ)$ is nontrivial. On the other hand, 
using Howald's characterization of multiplier ideals of monomial ideals
\cite{How}, one sees that
$$
\J(Y_e,cZ|_{Y_e}) = \J(Y_e,c_eZ|_{Y_e}) = (x_{e+1},\dots,x_n)
\fevery 1 \le e \le n-1,
$$
since, in view of \eqref{eq:eg:1}, for each fixed $e \ge 1$
we can neglect the term
$\frac{1}{a_1} + \dots + \frac{1}{a_e}$
in the computation of $\J(Y_e,cZ|_{Y_e})$. Therefore for every $e \ge 1$
the ideal $\J(Y_e,cZ|_{Y_e})$ is just the restriction
of the previous $\J(Y_{e-1},cZ|_{Y_{e-1}})$;
in other words, the multiplier ideals do not get any deeper.
Of course the situation is very different if the $Y_e$ are not general:
if for instance we take $Y_e = \{ x_{n-e+1} = \dots = x_n = 0\}$,
then we find that
$$
\J(Y_e,cZ|_{Y_e}) \subseteq (x_{e+1},\dots,x_n)^e
$$
which shows that for special restrictions
the multiplier ideals do indeed get deeper.
\end{eg}

This example gives an instance where inversion of adjunction fails
for non-effective pairs. It tells us that if we are looking for any reasonable
generalization of Theorem~\ref{thm:usual-adj}, then we must impose some
conditions on $Y$. 
As the following result shows, in the case of a homogeneous valuation
it suffices to impose a suitable tangency condition to $Y$
to gain control on the situation.

\begin{thm}\label{thm:discr-restr:homog-case}
Let $E$ be a prime divisor over $X = \A^n$ with center a point $P \in X$, 
and suppose that the valuation $\val_E$ is homogeneous with respect to 
some linear coordinates $\xx= (x_1,\dots,x_n)$ centered at $P$.
Let $Y = \A^{n-e} \subseteq X$ be a linear subspace of codimension $e < n$ that is
tangent to a principal tangent direction of $E$ at $P$
(cf. Definition~\ref{defi:principal-tangents}).
Then there exist a prime divisor $F$ over $Y$ with center $P$,
and an integer $q \ge 1$, such that
for every proper closed subscheme $Z \subset X$ not containing $Y$ and every $c \ge 0$
$$
q\,a_F(Y,cZ|_Y - eP) \le a_E(X,cZ).
$$
In particular, if $a_E(X,cZ) \ge 0$, then $a_F(Y,cZ|_Y - eP) \le a_E(X,cZ)$.
\end{thm}

\begin{proof}
Let
$$
W = W^1(E) \subset J_\infty X
$$
be the irreducible contact locus associated with $\val_E$, and let $\m = \val_E(P)$.
Our assumption is that $\T_PY$ contains a principal tangent 
vector $\x$ of $E$ at $P$.
By definition, $\x$ is the image in $T_PX$
of a non-constant jet $\z \in W_\m$ via the map 
$\ff_{\m,1}\colon T_P^{(\m)}X \to T_PX$. 
Since $\T_PY \cap \ff_{\m,1}(W_\m) = \ff_{\m,1}(J_\m Y \cap W_\m)$
by Lemma~\ref{lem:T_P^mX-ff-Y_m}, we can choose 
$$
\z \in J_\m Y \cap W_\m.
$$

The invariance of $\val_E$ and $W$ under the canonical respectively 
homogeneous $\C^*$-action imply that with $\g(t) \in W$
we also have $\l^{-\m}\. \g(\l t) \in W$.
In other words, the set $W$ is invariant under the action 
$\g(t) \mapsto \l^{-\m}\. \g(\l t)$. 
In the coordinates $x_i^{(p)}$, this action is given by
$x_i^{(p)} \mapsto \l^{p-\m}x_i^{(p)}$. 
Suppose that $\g \in W$ is an element mapping to $\z$. 
Wiring
$$
\g(t) =
\Big(\sum_{p \ge \m} a_1^{(p)} t^p,\dots,\sum_{p \ge \m} a_n^{(p)} t^p\Big),
$$
this means that $\z(t) = \big(a_1^{(\m)} t^\m,\dots,a_n^{(\m)} t^\m\big)$. 
Taking the limit for $\l \to 0$ under the above action, $\g$ degenerates to the arc
$\g^0(t) = \big(a_1^{(\m)} t^\m,\dots,a_n^{(\m)} t^\m\big)$,
which still belongs to $W$ since the latter is closed,
and maps to $\z$. Note that $\g^0$ also belongs to $J_\infty Y$. 
This proves that
$$
\z \in (J_\infty Y \cap W)_\m.
$$
We can thus select an irreducible component $C$ of $J_\infty Y \cap W$
whose image $C_\m$ in $J_\m Y$ contains $\z$.

Note that $C$ is an irreducible closed cylinder in $J_\infty Y$
that does not dominate $Y$, and
$$
P_\m \ne C_\m \subseteq \T_P^{(\m)}Y.
$$
By Theorem~\ref{thm:ELM}, there exists a
prime divisor $F$ over $Y$ and a positive integer $q$
such that $\val_C = q\val_F$,
as valuations of $\C(Y)$. The divisor $F$ has center $P$ in $Y$,
and
$$
q\val_F(P) = \m = \val_E(P).
$$
We have
$\min \{\ord_\g^X(Z) \mid \g \in W \} \le \min \{\ord_\g^X(Z) \mid \g \in C \}$
by the inclusion $C \subseteq W$, and $\ord_\g^X(Z) = \ord_\g^Y(Z|_Y)$
for every $\g \in C$ by Remark~\ref{rmk:order-restriction},
since $C \subset J_\infty Y$. Therefore 
$$
q\val_F(Z|_Y) =\val_C(Z|_Y) \ge \val_W(Z) = \val_E(Z).
$$
We now estimate the codimension of $C$ in $J_\infty Y$.
We can assume that $Y$ is defined by $x_i = 0$ for $0 \le i \le e$.
Since $W_{\m-1} = P_{\m-1}$, the functions
$x_i^{(p)}$, for $1\le i \le e$ and $\ 0 \le p \le \m - 1$,
vanish identically on $W$. 
It follows that $\codim(C,Y_\infty) \le \codim(W,X_\infty) - e\m$.
By Theorem~\ref{thm:ELM}, this implies that
$$
q\big(k_F(Y) + 1\big) \le k_E(X) + 1 - e\val_E(P).
$$
Altogether, we have
\begin{align*}
q\,a_F(Y,cZ|_Y - eP)
& = q\big(k_F(Y) + 1 - c\val_F(Z|_Y) + e\val_F(P)\big) \\
& \le k_E(X) + 1 - c\val_E(Z)  \\
& = a_E(X,cZ). \qedhere
\end{align*}
\end{proof}

In the language of multiplier ideals, the theorem implies
the following property.

\begin{cor}\label{cor:discr-restr:homog-case}
With the same notation and assumptions of Theorem~\ref{thm:discr-restr:homog-case},
suppose that $a_E(X,cZ) \le a$ for some $a \in \{0,1,\dots,e\}$.
Then 
$$
\mm^{e-a} \not\subseteq \J(Y,cZ|_Y).
$$
\end{cor}

\begin{proof}
Theorem~\ref{thm:discr-restr:homog-case} implies that
$a_F(Y,cZ|_Y-eP) \le a$ for some
prime divisor $F$ over $Y$ with center $P$, and therefore
$a_F(Y,cZ|_Y-(e-a)P) \le 0$ since $\val_F(P) \ge 1$.
This implies the formula in the statement.
\end{proof}

Theorem~B stated in the introduction follows as 
a particular case of this corollary, by taking $a=0$.

\section{Birational rigidity}\label{sec:rigidity}

The minimal model program predicts that
every variety of Kodaira dimension $-\infty$
is birationally equivalent to a Mori fiber space, namely,
a normal projective variety with terminal
$\Q$-factorial singularities, endowed with an
extremal Mori contraction of fiber type. 
The main motivation of the program is to choose this variety as the ``simplest'' model
to represent the birational equivalence class of the originally given variety.
It is then natural to inquire about the unicity of such a choice.

This led to the notions of birational rigidity and birational superrigidity. 
There are various versions of these notions in the literature. 
Here we adopt the definition of birational superrigidity given by Corti \cite{Co2} 
(which we recall next) as this is the one fitting most naturally
within the framework of the minimal model program. 
The next theorem also holds, however, 
with definition of birational superrigidity given in \cite{Pu2}
(see Remark~\ref{rem:defi-rigidity} below), and thus answers
the conjecture from \cite{Pu2}.

\begin{defi}
A Mori fiber space $X$ is said to be {\it birationally superrigid}
if every birational map from $X$ to any Mori fiber space
is an isomorphism preserving the fibration.
\end{defi}

This is a very strong condition. It implies that $X$ is the
only Mori fiber space in its birational equivalence class
and, in particular, that $X$ is not rational.
It also implies that the group of birational transformations
$\Bir(X)$ of $X$ coincides with the group of automorphisms $\Aut(X)$.
The slightly weaker notion of birational rigidity only requires the
Mori fiber structure to be ``birationally"
defined in a unique way; we refer the interested reader to \cite{Co2} for the
precise definition.

The following theorem, stated in the introduction as Theorem~A, 
classifies all birationally superrigid smooth hypersurfaces in projective spaces. 

\begin{thm}\label{thm:bir-rigidity}
For $N \ge 4$, every smooth hypersurface
$X \subset \P^N$ of degree $N$ is birationally superrigid.
\end{thm}

For the reader familiar with questions related to birational rigidity,
the basic idea of the proof is the usual one.
One starts by supposing that a birational map
$$
\f \colon X \rat X'
$$
between $X$ and some other Mori fiber space $X'$
is given, and assume for contradiction that $\f$ is not an isomorphism.
Generally speaking, the {\it method of maximal singularities} consists
in quantifying the singularities of
the indeterminacy of the map, aiming for
a contradiction in view of the (supposedly, relatively low) degrees of the variety
and of the equations defining the map.

The first key ingredient, known as the
{\it Noether--Fano inequality}, gives a precise bound on the singularities
of the indeterminacies of $\f$. In the case at hand
(i.e., when $X$ is a smooth hypersurface of degree $N$ in
$\P^N$, with $N \ge 4$), we can state this bound as follows.
Fix a projective embedding $X' \subseteq \P^m$, and let
$\H = \f^{-1}_*|\O_{X'}(1)|$ (in other words,
$\H$ is the linear system giving the map $X \rat \P^m$).
By the Lefschetz hyperplane theorem, we have
$$
\H \subseteq |\O_X(r)| \ \ \text{ for some } \ \ r \ge 1.
$$
Let $B(\H) \subset X$ be the base scheme of $\H$.

\begin{prop}[Noether--Fano Inequality \protect{\cite{IM}}]\label{NFI}
If $\f$ is not an isomorphism,
then the canonical threshold of the pair $(X,B(\H))$
satisfies the inequality
$$
\can{X,B(\H)} < 1/r.
$$
\end{prop}

This result is essentially due to \cite{IM}.
The general version of this property, holding when $X$ is an arbitrary
Mori fiber space, was given in dimension three in \cite{Co1}
and then extended in all dimensions in \cite{Is,Ma,dF}.
The proof of Theorem~\ref{thm:bir-rigidity} 
will be the result of the combination of the Noether--Fano inequality
with the following lower-bounds on canonical thresholds on hypersurfaces.

\begin{thm}\label{thm:bound-codim-2-case}
Let $V \subset \P^{n+1}$ be a smooth hypersurface of dimension $n \ge 3$
and degree $d$. Let $B \subset V$ be a proper closed subscheme of
codimension $\ge 2$, and assume that for some $r\ge 1$ the sheaf
$\O_V(r)\otimes I_B$ is globally generated. Then
\begin{align*}
\can{V,B} &\ge \min \bigg\{ \frac1r , \, \frac 2r\sqrt{\frac{n-2}{d}} \bigg\}
&\text{if $3 \le n \le 5$,} \\
\can{V,B} &\ge \min \bigg\{ \frac1r , \, \frac 1r\sqrt{\frac{2(n-2)}{d}} \bigg\}
&\text{if $n \ge 6$.}
\end{align*}
\end{thm}

\begin{rmk}
We expect the first bound to hold in general, for every $n \ge 3$. 
\end{rmk}

We will prove this theorem in the following sections.
Granting the theorem for the moment, we can quickly
finish the proof of the birational rigidity of $X$.

\begin{proof}[Proof of Theorem~\ref{thm:bir-rigidity}]
Keeping the notation already introduced, we observe that
$B(\H)$ has codimension $\geq 2$ in $X$ and
$\O_X(r)\otimes I_{B(\H)}$ is globally generated.
Combining Proposition~\ref{NFI}
with the second formula in Theorem~\ref{thm:bound-codim-2-case}, we obtain
$$
\frac{1}{r}\sqrt{\frac{2(N-3)}{N}} < \frac{1}{r}.
$$
This gives $N < 6$. We can therefore apply 
the first formula in Theorem~\ref{thm:bound-codim-2-case}, 
which implies that $N < 4$.
\end{proof}

\begin{rmk}\label{rem:defi-rigidity}
Theorem~\ref{thm:bir-rigidity} holds no matter which definition of birational
rigidity is adopted. The proof goes indeed by showing that, 
in the language of \cite{Pu2}, there are no {\it maximal singularities} on $X$, and
the same property also implies birational superrigidity in the sense 
defined in \cite{Pu2}. 
Specifically, we prove that for any movable linear system $\H \subset |\O_X(r)|$
the canonical threshold of the pair $(X,B(\H))$ is bounded below by $1/r$. 
Now, $r$ is equal to the {\it threshold of canonical adjunction}
of $(X,\H)$ defined in \cite{Pu2}, and
\cite[Section~2, Proposition~1]{Pu2}, the version of Noether--Fano
inequality needed there, holds for all smooth hypersurfaces 
of degree $N$ in $\P^N$, with $N \ge 4$ (the proof of this proposition 
does not use the hypothesis of $(N-1)$-regularity: such condition is only needed
later, to exclude the existence of maximal singularities).
In particular, Theorem~\ref{thm:bir-rigidity}
answers positively the main conjecture of \cite{Pu2} and implies its main result
as a special case. 
\end{rmk}

\section{Few more preliminaries for the proof of Theorem~\ref{thm:bound-codim-2-case}}

We gather in this section a few more properties that we will need in
the proof of Theorem~\ref{thm:bound-codim-2-case}.
The first property concerns the behavior of log canonical thresholds via
generic projections. The statement is a slightly more precise 
formulation of \cite[Theorem~1.1]{dFEM2}, which follows entirely
from the original proof. We should however point out a difference of notation: 
here the symbol $a_E(\;\,)$ denotes log-discrepancies, whereas
in \cite{dFEM2} the same symbol was used to denote discrepancies.

\begin{thm}[\protect\cite{dFEM2}]\label{thm:projection}
Let $E$ be a divisor over a smooth variety $X$, with center a point $P \in X$,
and let $Z \subseteq X$ be a Cohen-Macaulay subscheme of codimension $k \ge 1$.
Assume that $f \colon X \to Y$ is a smooth and
proper morphism onto a smooth variety $Y$
of dimension $\dim Y = \dim Z + 1$ such that $f|_Z$ is a finite map.
Then the valuation $\val_E$ restricts to a divisorial valuation
$q \val_F$ on $\C(Y)^*$ where $F$ is a divisor over $Y$ and $q \ge 1$, 
and
$$
q\,a_F\big(Y,b(c/k)^kf_*[Z]\big) \le a_E(X,cZ)
$$
where $b = 1$ if $Z$ is locally complete intersection and $b = k!$ otherwise.
\end{thm}

Notice that $X$ needs not be proper or projective in this theorem.

\begin{rmk}\label{rmk:projection}
Since $\val_F$ is obtained from the restriction of $\val_E$, it follows that the
center of $F$ on $Y$ is $Q :=f(P)$. Moreover, if
$\^\G_E$ is not in the kernel of the map
$df|_P \colon T_PX \to T_QY$, then 
$$
\^\G_F = df|_P(\^\G_E) \subseteq \T_QY
$$ 
(cf. Notation~\ref{not:G_E}). Finally, if $f$ is a linear map of affine spaces
and $\val_E$ is an homogeneous valuation, then $\val_F$ is homogeneous too. 
\end{rmk}

Next we review a few properties on multiplicities from \cite{dFEM2}.
We begin with some basic facts from intersection theory.
The {\it multiplicity} of a pure-dimensional scheme $Z$ at a point $P$
is defined to be $e_P(Z) := e(\mm_{Z,P})$, the Samuel multiplicity of
the maximal ideal $\mm_{Z,P} \subset \O_{Z,P}$, as in \cite[Example~4.3.4]{Fu}. 
If $\a = \sum n_i[V_i]$ is a cycle on a variety $X$, where each
$V_i$ is a subvariety, then 
we define the multiplicity of $\a$ at a point $P \in X$
to be $e_P(\a) := \sum n_i \,e_P(V_i)$, where
we set $e_P(V_i) = 0$ if $P \not\in V_i$.

\begin{rmk}\label{rmk:fund-cycle}
These definitions are compatible:
if $Z$ is a pure-dimensional closed subscheme of a variety $X$, and $[Z]$ is the
associated fundamental cycle, then $e_P(Z) = e_P([Z])$ for every point
$P \in Z$ (see \cite[Example~4.3.4]{Fu}).
\end{rmk}

\begin{rmk}\label{rem:7.4}
If $D$ is an effective Cartier divisor on a variety $X$ and $P \in X$ is a regular point, 
then $e_P(D)$ is simply the multiplicity of a
generator of the ideal of $D$ in the local ring at $P$
(see \cite[Example~4.3.9]{Fu}).
\end{rmk}

Let $Z$ be a pure-dimensional scheme. The function
$P \mapsto e_P(Z)$ is upper-semicontinuous by \cite[Theorem~(4)]{Ben}. 
In particular, if $T \subset Z$ is a irreducible closed set, 
then we can define the multiplicity of $Z$ along $T$ to be
\[
e_T(Z) := \min_{P \in T} e_P(Z),
\]
and we have $e_T(Z) = e_P(Z)$ for a general $P \in T$. 
Again, this notion extends to cycles:
if $\a = \sum n_i[V_i]$ is a cycle on a variety $X$, and $T \subseteq X$
is any irreducible closed subset, then we define $e_T(\a) := \sum n_i \,e_T(V_i)$, where
we set $e_T(V_i) = 0$ if $T \not\subseteq V_i$.
Note that $e_T(\a) = \min_{P \in T} e_P(\a)$ and the minimum
is achieved for general $P \in T$.

\begin{prop}\label{prop:restr-multiplicities}
Let $Z$ be  closed Cohen-Macaulay subscheme of $\P^m$ of
positive pure-dimension, and let $H \subseteq \P^m$ be a hyperplane.
\begin{enumerate}
\item[(i)]
If $H$ meets properly the embedded tangent cone
of $Z$ at a point $P$, then $e_P(Z \cap H) = e_P(Z)$.
\item[(ii)]
For every given hyperplane $\H \subset (\P^m)^\vee$,
if $H \in \H$ is general enough, then $e_P(Z \cap H) = e_P(Z)$
for every $P \in Z \cap H$.
\end{enumerate}
\end{prop}

\begin{proof}
We can assume that $Z \ne \P^m$.
Consider any linear subspace $L \subset \P^m$ of dimension
$\dim L = m - \dim Z$ that meets properly the embedded tangent cone
of $Z$ at $P$. Then the component of $Z \cap L$ at $P$
is zero-dimensional, and we have $e_P(Z) = \length(\O_{P,Z\cap L})$
by \cite[Proposition~7.1 and Corollary~12.4]{Fu}. This implies (i).
As long as $H$ is chosen outside the dual variety $Z_i^\vee$
of each irreducible component $Z_i$ of $Z$, (ii) follows from (i),
as it is explained for instance in \cite[Proposition~4.5]{dFEM2}.
Therefore, to conclude, it suffices to observe that $Z_i^\vee$ cannot contain
any hyperplane of $(\P^m)^\vee$, since it is irreducible of dimension $\le m-1$, and
$Z_i^{\vee\vee} = Z_i$ is not a point.
\end{proof}

\begin{prop}\label{prop:proj-multiplicities}
Let $Z$ be a pure-dimensional closed subscheme of $\P^m$.
Let $\f \colon \P^m\setminus\L \to \P^k$ be
a linear projection from a center $\L$ disjoint from $Z$, and assume that
$\f|_{Z_{\rm red}}$ is injective over the image of a point $P \in Z$.
If $\f^{-1}(\f(P))$ meets properly the embedded tangent cone of $Z$ at $P$,
then $e_P(Z) = e_{\f(P)}(\f_*[Z])$.
\end{prop}

\begin{proof}
By Remark~\ref{rmk:fund-cycle}, we can reduce
to the case in which $Z$ is a subvariety
of $\P^m$. Then the assertion is proven in the second part
of the proof of \cite[Proposition~4.6]{dFEM2}.
\end{proof}

\begin{prop}[\protect\cite{Pu3}]\label{prop:Puk}
Let $X \subset \P^{n+1}$ be a hypersurface of dimension $n \ge 2$, and let $\a$ be an
effective cycle on $X$ of pure codimension $k$, with $1 \le k \le n/2$, such that
$\a \equiv m\.c_1(\O_X(1))^k$ for some $m \in \N$.
Then $e_T(\a) \le m$ for every closed subvariety $T \subseteq X$ of dimension
$\ge k$ and not meeting the singular locus of $X$.
\end{prop}

\begin{proof}
For $k < n/2$ this is \cite[Proposition~5]{Pu3},
and the property extends to the extremal case $k = n/2$
(cf. \cite[Remark~4.4]{dFEM2}).
\end{proof}

Finally, we recall the following well-known fact, whose proof is here
recalled for the benefit of the reader.

\begin{prop}\label{prop:mult-canonical}
Let $A$ be an effective $\Q$-divisor on a smooth variety $X$,
and suppose that $a_E(X,A) \le 1$ for some prime exceptional divisor $E$
over $X$. If $T$ is the center of $E$ in $X$, then $e_T(A) \ge 1$.
\end{prop}

\begin{proof}
We can assume that $E$ is an exceptional divisor
of a log-resolution $f \colon X' \to X$ of $(X,A)$.
Pick a general point $P \in T$, and let $Y \subset X$ be a general
complete intersection subvariety of codimension $\codim(Y,X) = \dim T$,
passing through $P$. Then
the proper transform $Y'$ of $Y$ meets $E$ transversally, and we have
$a_{E'}(Y,A|_Y) \le 1$ if $E'$ is a component of $E|_{Y'}$.
Notice that $\dim Y \ge 2$.
If $H \subset Y$ is a general hyperplane section through $P$,
then $(H,A|_H)$ is not log-terminal at $P$ by Theorem~\ref{thm:usual-adj}.
Taking a general (smooth) complete intersection curve $C \subset H$ through $P$, we see
that $(C,A|_C)$ is not 
Kawamata log-terminal at $P$ by the same theorem. This is equivalent to
$e_P(A|_C) \ge 1$. On the other hand, by taking the hyperplanes cutting out $C$
generally enough, we can ensure that $e_P(A|_C) = e_P(A)$
(see Remark~\ref{rem:7.4} or Proposition~\ref{prop:restr-multiplicities}(i)).
We conclude that $e_T(A) \ge 1$.
\end{proof}

\section{Proof of Theorem~\ref{thm:bound-codim-2-case}}

Suppose that $V \subset \P^{n+1}$ is a smooth hypersurface of dimension $n \ge 3$
and degree $d$, $B \subset V$ is a proper closed subscheme of
codimension $\ge 2$, and the sheaf
$\O_V(r)\otimes I_B$ is globally generated for some $r\ge 1$.

After replacing $B$ with the intersection of two general
members of $|\O_V(r)\otimes I_B|$, we reduce to prove the
theorem when $B$ is a complete intersection subscheme
of $V$ cut out by two forms of degree $r$. Suppose that
$$
c := \can{V,B} < 1/r.
$$
For every $D \in |\O_V(r)\otimes I_B|$ the pair
$(V,cD)$ is terminal in dimension one 
by Propositions~\ref{prop:Puk} and~\ref{prop:mult-canonical}.
This implies that $(V,cB)$ is also terminal
in dimension one. Therefore the canonical threshold of
$(V,B)$ is computed by some divisor over $V$ with center equal
to a point $P \in V$.
Let $V' \subset V$ be a general hyperplane section through $P$, and let
$B' := B|_{V'}$. Note that
$$
V' \subset \P^n
$$
is a smooth hypersurface of degree $d$. By Theorem~\ref{thm:usual-adj},
there exists a prime divisor $E$ over $X$ with center $P$ such that
$$
a_E(V',cB') \le 0.
$$
If the choice of the hyperplane section
cutting $V'$ is sufficiently general, then we have $e_T(B') \le r^2$
for any positive dimensional closed subvariety $T \subseteq V'$
by Propositions~\ref{prop:Puk} and~\ref{prop:restr-multiplicities}(ii).

We consider a general linear projection 
$$
\r \colon \P^n \rat \P^{n-2}.
$$ 
Let $f \colon V' \rat \P^{n-2}$ be the restriction of $\r$ to $V'$.
We assume that the center of projection is disjoint from $B'$.
Let $A := f_*[B']$ and $Q := f(P)$.
Note that $A$ is an effective divisor of degree $dr^2$ on $\P^{n-2}$.
By Theorem~\ref{thm:projection}, there is a divisor $F$ over $\P^{n-2}$
with center $Q$ and such that
$$
a_F(\P^{n-2},\l A) \le 0  \for  \l = \frac{c^2}{4}.
$$

If $3 \le n \le 5$, then one can argue as in \cite{dFEM2} to see that, near $Q$, 
the pair $(\P^{n-2},\l A)$ is Kawamata log-terminal 
away from $Q$. One can then conclude quickly as follows. 
Nadel's vanishing theorem \cite[Theorem~9.4.8]{Laz} gives
$$
\HH{1}{\P^{n-2},\O_{\P^{n-2}}(\lru \l dr^2 + \ep \rru - n+1)\otimes\J(\P^{n-2},\l A)} = 0
\fall \ep > 0,
$$
This implied that, if $\S$ is the irreducible component supported at $Q$
of the scheme defined by the multiplier ideal $\MI{\P^{n-2},\l A}$, then there is a surjection
$$
\HH{0}{\P^{n-2},\O_{\P^{n-2}}(\lru \l dr^2 + \ep \rru - n+1)} \surj \HH{0}{\P^{n-2},\O_\S}
\fall \ep > 0.
$$
It follows that
$\lru \l dr^2 + \ep \rru - n+1 \ge 0$ for all $\ep > 0$.
Recalling that $\l = c^2/4$, we get
$$
c \ge \frac 2r\sqrt{\frac{n-2}{d}}.
$$
This gives the first inequality in the statement of the theorem.

The problem when $n$ gets larger is that one looses control on the multiplicities
of $A$ (cf. \cite[Remark~4.2]{dFEM2}). 
To deal with this problem, we restrict to the affine setting
and consider flat degenerations to homogeneous schemes in order
to reduce to a setting where we can apply
Theorem~\ref{thm:discr-restr:homog-case} to cut down the dimension. 
We will perform the projection in two steps, by breaking it
into the composition of two general linear projections. 
The degeneration to homogeneous schemes will be taken half way down, 
after the first projection. 
Theorem~\ref{thm:discr-restr:homog-case}
will then be applied after the second projection.
We proceed in several steps. 

\begin{step}
Let $H \subset \P^n$ be a general hyperplane, and let $\A^n = \P^n \smallsetminus H$. 
We assume in particular that $H$ intersect property every irreducible component of $B'$.
We denote by $X$ and $Z$ the affine traces of $V'$ and $B'$ in $\A^n$, 
so that 
$$
X \subset \A^n
$$
is a smooth hypersurface of degree $d$, $Z$ is a closed subscheme
of codimension two of $X$ cut out by two polynomials of degree $r$.
Note that $B'$ is equal to the projective closure $\ov Z$ of $Z$ in $\P^n$.

Recall that 
$$
a_E(X,cZ) \le 0
$$
for some divisor $E$ with center $P$ in $X$. Let
$$
W = W^1(E) \subset J_\infty X
$$
be the contact locus of $E$. We henceforth fix an integer $m$ such
that $W$ is a cylinder over $W_m$. 
By viewing $W_m$ as a closed subset of $J_m\A^n$, we can take its
homogeneous limit 
$$
(W_m)^0 \subset J_m\A^n
$$ 
with respect to the homogeneous $\C^*$-action on $\A^n$ 
along some linear coordinates $\xx = (x_1,\dots,x_n)$
centered at $P$, as described in  
Section~4.
Note that $(W_m)^0$ is contained in the fiber over $P$. 
Let
$$
\m := \min \{\, q \ge 0 \mid ((W_m)^0)_q \ne P_q \,\}.
$$
Then $((W_m)^0)_\m \subset T_P^{(\m)}\A^n$ and,
recalling the isomorphism $\ff_{\m,1} \colon T_P^{(\m)}\A^n \to T_P\A^n$, 
we can fix a nonzero tangent vector
$$
\x \in \ff_{\m,1}\big(((W_m)^0)_\m\big) \subset T_P\A^n.
$$
\end{step}

\begin{step}
We now consider a general linear projection
$$
\s \colon \A^n \to \A^{n-1}.
$$
We assume that
the induced map $g\colon X \to \A^{n-1}$ is \'etale in a neighborhood of $P$. 
The valuation $\val_E$ on $X$ restricts to a valuation $\val_{E'}$
on $\A^{n-1}$ where $E'$ is a divisor with center $P' = g(P)$.
The map $g_\infty \colon J_\infty X \to J_\infty\A^{n-1}$ induces
an isomorphism from $\p^{-1}(P)$ to $(\p')^{-1}(P')$,
where $\p \colon J_\infty X\to X$ and $\p' \colon J_\infty \A^{n-1}\to \A^{n-1}$ 
are the two projections. Under this isomorphism,
$W$ maps to the contact locus
$$
W' := W^1(E')\subset J_\infty \A^{n-1}.
$$
We consider the homogeneous limit
$$
(W')^0 \subset J_\infty \A^{n-1}
$$
of $W$ with respect to the homogeneous 
$\C^*$-action of $\A^{n-1}$ on some linear coordinates $\yy = (y_1,\dots,y_{n-1})$
centered at $P'$. Note that $\s$ is compatible with the two $\C^*$-actions.

\begin{lem}\label{lem:Winfty}
$\s_\m\big(((W_m)^0)_\m\big) = ((W')^0)_\m$.
\end{lem}

\begin{proof}
Since the maps $\s_q$ commute with the various projections
on the jet schemes of $\A^n$ and $\A^{n-1}$,
it is enough to prove that $\s_m\big((W_m)^0\big) = ((W')^0)_m$.
In fact, since $W'$ is the isomorphic image of $W$, 
$W'$ is a cylinder over $W'_m$ and hence
$((W')^0)_m = (W'_m)^0$ by Lemma~\ref{lem:W}.
Therefore, all we need to prove is that
\begin{equation}\label{eq:32}
\s_m\big((W_m)^0\big) = (W'_m)^0. 
\end{equation}

Let $\p_{m,0} \colon J_mX \to X$, 
$\~\p_{m,0} \colon J_m\A^n \to \A^n$ and
$\p'_{m,0} \colon J_m \A^{n-1} \to \A^{n-1}$
denote the various projections.
The map $\s_m$ restricts to a linear map
$$
\f \colon \~\p_{m,0}^{-1}(P) = \Spec\C[\xx',\dots,\xx^{(m)}] \to
(\p_{m,0}')^{-1}(P') = \Spec\C[\yy',\dots,\yy^{(m)}].
$$
We consider the induced rational map on projective spaces
$$
\ov{\f} \colon \Proj\C[u,\xx',\dots,\xx^{(m)}] \rat
\Proj\C[u,\yy',\dots,\yy^{(m)}].
$$
For any subset $S \subseteq \~\p_{m,0}^{-1}(P)$, we denote by
$\ov{S}$ its closure in $\Proj\C[u,\xx',\dots,\xx^{(m)}]$.
Note that $\p_{m,0}^{-1}(P) \subset \~\p_{m,0}^{-1}(P)$, so we can consider 
its closure $\ov{\p_{m,0}^{-1}(P)}$ in $\Proj\C[u,\xx',\dots,\xx^{(m)}]$. 

We can assume that $\{x_n = 0\}$ is tangent to $X$ at $P$
and $\s$ is given by $y_i = x_i$ for $i=1,\dots,n-1$, so that 
the base locus $\L_m \subset \Proj\C[u,\xx',\dots,\xx^{(m)}]$ 
of $\ov{\f}$ is defined by $u=0$ and $x_i^{(p)} = 0$ for
all $p$ and all $i \ne n$.
By Proposition~\ref{prop:cone-proj}, the equation
\eqref{eq:32} holds if
$\ov{W_m} \cap \L_m = \emptyset$.
Since $W_m \subseteq \p_{m,0}^{-1}(P)$, it suffices to show that
\begin{equation}\label{eq:Winfty}
\ov{\p_{m,0}^{-1}(P)} \cap \L_m = \emptyset.
\end{equation}

We prove \eqref{eq:Winfty} using induction on $m$. The statement is trivial when $m=0$, 
so we can assume that $m \ge 1$ and $\ov{\p_{m-1,0}^{-1}(P)} \cap \L_{m-1} = \emptyset$. 
Proceeding by contradiction, suppose that \eqref{eq:Winfty}
is false for some $m$, and pick 
$$
\z \in \ov{\p_{m,0}^{-1}(P)} \cap \L_m.
$$
Let $(a_u, \underline{a}', \dots,\underline{a}^{(m)})$ denote the coordinates
of $\z$ in the coordinate system $(u,\xx',\dots,\xx^{(m)})$. 
Since $\z \in \L_m$, we have $a_u = 0$ and $a_i^{(p)} = 0$ for all $p$
and all $i \ne n$. Consider the rational map
$$
\ov{\om} \colon \Proj\C[u,\xx',\dots,\xx^{(m)}]
\rat \Proj\C[u,\xx',\dots,\xx^{(m-1)}].
$$
induced by the projection $\om \colon \~\p_{m,0}^{-1}(P) \to \~\p_{m-1,0}^{-1}(P)$. 
Let $\S$ denote the indeterminacy locus of $\ov{\om}$.

Note that $\ov{\om}$ maps $\L_m \smallsetminus \S$ to $\L_{m-1}$
and $\ov{\p_{m,0}^{-1}(P)} \smallsetminus \S$ to the closure $\ov{\p_{m-1,0}^{-1}(P)}$
of $\p_{m-1,0}^{-1}(P)$ in $\Proj\C[u,\xx',\dots,\xx^{(m-1)}]$.
Since \eqref{eq:Winfty} holds for $m-1$, it follows that $\z$
must lie in $\S$. This means that $a_n^{(p)} = 0$ of all $p < m$, and hence
the only nonzero coordinate of $\z$ is $a_n^{(m)}$. In particular, 
$$
\z = (0,\dots,0,a_n^{(m)}) \in \ov{T_P^{(m)}\A^n}.
$$
On the other hand, we have
$$
\ov{\p^{-1}_{m,0}(P)} \cap \ov{T_P^{(m)}\A^n} = \ov{T_P^{(m)}X} 
= \{x_n^{(m)}=0\} \subset \ov{T_P^{(m)}\A^n}.
$$
This gives a contradiction.
\end{proof}

If $\s$ is sufficiently general, then $((W_m)^0)_\m$ is not contained
in the kernel of the linear map $T_P^{(\m)}\A^n \to T_{P'}^{(\m)}\A^{n-1}$, 
and therefore Lemma~\ref{lem:Winfty} implies that
$$
\m = \m_{P'}\big((W')^0\big).
$$
Moreover, since the maps induced by $\s$ on the jet schemes commute with the 
isomorphisms $\ff_{m,1} \colon T_P^{(m)}\A^n \to T_P\A^n$ and
$\ff'_{m,1} \colon T_{P'}^{(m)}\A^{n-1} \to T_{P'}\A^{n-1}$,
we also see that the image of $\ff_{\m,1}\big(((W_m)^0)_\m\big)$
under the map
$$
d\s|_P \colon T_P\A^n \to T_{P'}\A^{n-1}
$$
is equal to $\ff'_{\m,1}\big(((W')^0)_\m\big)$.
We previously fixed a nonzero tangent vector
$$
\x \in \ff_{\m,1}\big(((W_m)^0)_\m\big) \subset T_P\A^n. 
$$
We consider its image 
$$
\x' \in \ff'_{\m,1}\big(((W')^0)_\m\big) \subset T_{P'}\A^{n-1}. 
$$
If the projection $\s$ is sufficiently general then $\x'$ 
will also be a nonzero vector. 

Let $C$ be an irreducible component of $(W')^0$ 
such that $\x' \in \ff'_{\m,1}(C_\m)$, and
let $\val_C$ be the associated valuation on $\A^{n-1}$.
By Theorem~\ref{thm:ELM}, we can find a divisor $F$ over $\A^{n-1}$
and an integer $q$ such that $\val_C = q\val_F$, and 
$\codim(W^q(F),J_\infty\A^{n-1}) \le \codim(C,J_\infty\A^{n-1})$.
Note that $F$ has center $P'$ in $\A^{n-1}$,
and $\x' \in \^\G_F \subseteq T_{P'}\A^{n-1}$
(cf. Notation~\ref{not:G_E}). 
Note also that $\val_F$ is a homogeneous valuation.

Since the center of the projection $\P^n \rat \P^{n-1}$ corresponding to $\s$ is
disjoint from the projective closure $\ov Z$ of $Z$, 
the induced morphism $Z \to \A^{n-1}$ is finite, and in particular the scheme image 
$$
Z' := \s(Z) \subset \A^{n-1}
$$ 
of $Z$ is a closed subscheme of $\A^{n-1}$. We consider its homogeneous limit
$$
(Z')^0 \subset \A^{n-1}
$$
under the $\C^*$-action above. Recall that $g \colon X \to \A^{n-1}$ is \'etale
in a neighborhood of $P$, and note that $Z' = g(Z)$. Then, 
by Proposition~\ref{prop:semicont}, we get
$$
a_F(\A^{n-1},c (Z')^0) \le \frac 1q a_{E'}(\A^{n-1},cZ') = \frac 1q a_E(X,cZ) \le 0.
$$

\begin{lem}
\label{lem:CM}
$(Z')^0$ is a Cohen--Macaulay scheme.
\end{lem}

\begin{proof}
First observe that since the projective closure $\ov Z$
of $Z$ is equal to the complete intersection subscheme $B' \subset \P^n$, 
the homogeneous limit $Z^0 \subset \A^n$ of $Z$ is the cone over $\ov Z \cap H$,
and is complete intersection in $\A^n$ (see Remark~\ref{rmk:Z^0-cone}).
Note also that $\ov Z$ is disjoint from
the center of the projection $\P^n \rat \P^{n-1}$, and thus
$Z^0 \to \A^{n-1}$ is a finite map with image
$$
\s(Z^0) = (Z')^0,
$$ 
by Proposition~\ref{prop:cone-proj}.
It is a general fact that the image under a finite morphism of a Cohen--Macaulay
scheme is Cohen--Macaulay: if $R \subset S$ is a module-finite extension of local rings 
then $\dim(R) = \dim (S)$ and 
$\depth(R) \ge \depth_R(S) = \depth(S)$ (see \cite[Proposition~III(3.16)]{AK}),
and thus $R$ is a Cohen--Macaulay whenever $S$ is Cohen--Macaulay.
Since $Z^0$ is Cohen--Macaulay, we conclude that $(Z')^0$ is Cohen--Macaulay too.
\end{proof}
\end{step}

\begin{step}
We next take a general linear projection 
$$
\t \colon \A^{n-1} \to \A^{n-2}.
$$
We consider the divisor on $\A^{n-2}$ given by
$$
A^0 := \t_*[(Z')^0].
$$
Note that $A^0 = \r_*[Z^0]$, where we now denote $\r = \t\o\s \colon \A^n \to \A^{n-2}$.
Let
$$
\x'' \in \T_Q\A^{n-2}
$$
the image of $\x'$ under the map 
$d\t|_{P'} \colon \T_{P'}\A^{n-1} \to \T_Q\A^{n-2}$, where $Q := \t(P')$.

Since $(Z')^0$ is Cohen--Macaulay, 
Theorem~\ref{thm:projection} implies that there exists
a divisor $G$ over $\A^{n-2}$ such that   
$$
a_G(\A^{n-2},\d A^0) \le 0 \for \d = \frac{c^2}{2}.
$$
As discussed in Remark~\ref{rmk:projection}, $G$ has center $Q$ in $\A^{n-2}$, 
$\val_G$ is homogeneous, and
$$
\^\G_G = d\t|_{P'}(\^\G_F) \subset T_Q\A^{n-2}. 
$$ 
In particular, $\x'' \in \^\G_G$.

\begin{lem}\label{lem:e_L<r^2}
If $L'' \subset \A^{n-2}$ is the line through $Q$
with tangent direction $\x''$, then 
$$
e_{L''}\big(A^0\big) \le r^2.
$$
\end{lem}

\begin{proof}
Let us denote by $\ov Z$ and $\ov{A^0}$ the closures of $Z$ and $A^0$
respectively in
$\P^n$ and $\P^{n-2}$. Note that $\ov Z = B'$. Also, let
$H'' = \P^{n-2} \smallsetminus \A^{n-2}$. 
Since $\r$ is the composition of general projections, the induced rational map
$\r_H\colon H \rat H''$ can be treated as a general projection.
The scheme $Z^0$ coincides with the affine cone over the intersection $\ov Z \cap H$, 
since $\ov Z$ is a complete intersection
subscheme of $\P^n$ and $H$ meets it properly (see Remark~\ref{rmk:Z^0-cone}). 
Similarly, $A^0$ is the affine cone over $\ov{A^0} \cap H''$. 
Note also that $(\r_H)_*[\ov Z \cap H] = \ov{A^0} \cap H''$.

Recall that the multiplicities of $\ov Z$ can exceed $r^2$ only over a finite set.
Since $H$ was chosen to be general, we can assume that 
$$
e_R(\ov Z\cap H) = e_R(\ov Z) \le r^2
$$ 
for every point $R \in \ov Z \cap H$,
see Proposition~\ref{prop:restr-multiplicities}(ii).

Suppose now that $R = L \cap H$ where $L \subset \A^n$ is the line through $P$
with tangent direction $\x$. 
Note that $R'' := \r_H(R)$
is the point of intersection of $L''$ with $H''$. 
If $R \not\in \ov Z \cap H$, then $\ov{A^0} \cap H''$ does not contain its image $R''$
for a general projection $\r|_H$, and so 
$$
e_{L''}\big(A^0\big) = e_{R''}(\ov{A^0}\cap H'') = 0.
$$
Suppose then that $R \in \ov Z \cap H$. A general plane  $M \subset H$ passing through $R$
will meet properly the tangent cone of $Z \cap H$ at $R$,
and will have no other intersections with $Z \cap H$ outside $R$
since $Z \cap H$ has codimension three in $H$. Therefore,
by assuming that the source of the projection $\r|_H$ is a general line in $H$,
we have 
$$
e_{L''}\big(A^0\big) = e_{R''}(\ov{A^0} \cap H'') = e_{R}(\ov Z \cap H) \le r^2
$$
by Proposition~\ref{prop:proj-multiplicities}.
\end{proof}
\end{step}

\begin{step}
By Proposition~\ref{lem:C_E-val_E(P)} (see also Definition~\ref{defi:principal-tangents}), 
we can find a principal tangent vector $\x^*$ of $G$ at $Q$ 
that is close enough to $\x''$ so that if $L^* \subset \A^{n-2}$ is the line
tangent to $\x^*$ at $Q$, then the bound
\begin{equation}\label{eq:e_L}
e_{L^*}(A^0) \le r^2. 
\end{equation}
still holds.

We take a general plane $\A^2 \subset \A^{n-2}$ containing $L^*$, 
and consider the $\Q$-divisor
$$
\D := \d A^0|_{\A^2}
$$
where $\d = c^2/2$ (as defined in Step~3).
In view of~\eqref{eq:e_L} and Proposition~\ref{prop:restr-multiplicities},
we have $\ord_C(\D) \le \d r^2$ for every irreducible curve $C \subset \A^2$
as long as the plane is sufficiently general.
Since $\d < 1/r^2$, this implies that the multiplier ideal $\J(\A^2,\D)$
defines a zero-dimensional subscheme of $\A^2$, which is nonempty at $Q$.
In fact, if $\S$ is the component of this scheme supported at $Q$,
$(x,y)$ are some affine coordinates of $\A^2$
centered at $Q$, and $I_\S \subset \C[x,y]$ is the ideal of $\S$ in $\A^2$, then
Corollary~\ref{cor:discr-restr:homog-case} implies that
$$
(x,y)^{n-4} \not\subseteq I_\S.
$$
Notice on the other hand that $\J(\A^2,\D)$, and so $I_\S$, are homogeneous
ideals in $\C[x,y]$. It follows that, among the generators
of $\O_\S = \C[x,y]/I_\S$, 
one needs some homogeneous forms in $\C[x,y]$ of degree at least $n-4$.

Consider the component-wise closure $\ov{\D}$ of $\D$ in $\P^2$. 
By construction, $\ov{\D}$ is an effective $\Q$-divisor
of degree $\d dr^2$. Nadel's vanishing theorem gives
$$
\HH{1}{\P^2,\O_{\P^2}(\lru \d dr^2 + \ep \rru - 3)\otimes\J(\P^2,\ov{\D})} = 0
\fall \ep > 0.
$$
Since $\S$ is a zero-dimensional connected component
of the scheme defined by $\J(\P^2,\ov{\D})$,
this yields the surjection
$$
\HH{0}{\P^2,\O_{\P^2}(\lru \d dr^2 + \ep \rru - 3)} \surj \HH{0}{\P^2,\O_\S}
\fall \ep > 0.
$$
By what just observed on the generators of $\O_\S$, we conclude that
$\lru \d dr^2 + \ep \rru - 3 \ge n-4$ for all $\ep > 0$.
Recalling that $\d = c^2/2$, we get
$$
c \ge \frac 1r\sqrt{\frac{2(n-2)}{d}}, 
$$
which gives the second inequality in the statement of the theorem.
This completes the proof of Theorem~\ref{thm:bound-codim-2-case}.
\end{step}

\providecommand{\bysame}{\leavevmode \hbox \o3em
{\hrulefill}\thinspace}

\end{document}